\documentclass{article}
 \usepackage{amstex}
 \usepackage{theorem}

 \newcommand\Oh{{\cal O}} % structure sheaf
 \newcommand\C{{\Bbb C}}
 \newcommand\Q{{\Bbb Q}}
 \newcommand\Z{{\Bbb Z}}
 \newcommand\R{\Bbb R}    % real numbers
 \newcommand\om{\omega}
 \newcommand\ep{\varepsilon}
 \newcommand\sI{{\cal I}}
 \newcommand\codim{\operatorname{codim}}
 \newcommand\Bl{\operatorname{Bl}}
 \newcommand\Coh{\operatorname{Coh}}
 \newcommand\Hom{\operatorname{Hom}}
 \newcommand\Hilb{\operatorname{Hilb}}
 \newcommand\GHilb{\operatorname{{\mbox{$G$}}-Hilb}}
 \newcommand\Lie{\operatorname{Lie}}
 \newcommand\Pic{\operatorname{Pic}}
 \newcommand\Stab{\operatorname{Stab}}
 \newcommand\card{\operatorname{\#}} % cardinality
 \newcommand\proj{{\Bbb P}}
 \newcommand\age{\operatorname{age}} % age
 \newcommand\im{\operatorname{im}} % image
 \newcommand\centre{\operatorname{centre}} % centre
 \newcommand\disc{\operatorname{disc}} % discrepancy
 \newcommand\into{\hookrightarrow}
 \newcommand{\broken}{\mathrel{{\relbar\kern-.2pt\rightarrow}}}
 \newcommand\tensor{\otimes}
 \newcommand\al{\alpha}
 \newcommand\be{\beta}
 \newcommand\ga{\gamma}
 \newcommand\de{\delta}
 \newcommand\la{\lambda}
 \newcommand\si{\sigma}
 \newcommand\De{\Delta}
 \newcommand\Si{\Sigma}
 \newcommand\Om{\Omega}

 \newcommand\sE{{\cal E}} % sheaf E
 \newcommand\sF{{\cal F}} % sheaf F
 \newcommand\sG{{\cal G}} % sheaf G
 \newcommand\sL{{\cal L}} % invertible sheaf L
 \newcommand\sZ{{\cal Z}} % universal family of clusters

 \newcommand{\Span}[1]{\left<#1\right>}        % <x> span or hull of x

 \newcommand\QED{\ifhmode\unskip\nobreak\fi\quad {\rm Q.E.D.}} % QED
 \newcommand\GL{\operatorname{GL}}      % full linear group
 \newcommand\SL{\operatorname{SL}}      % special linear group
 \newcommand\iso{\cong}
 \newcommand\1{^{-1}} % inverse
 \DeclareMathSymbol{\onto}     {\mathrel}{AMSa}{"10}   % -->>
 \newcommand\bij{\leftrightarrow}           % <-->
 \DeclareMathSymbol{\boxtimes}     {\mathbin}{AMSa}{"02}
 \newcommand\id{\mathrm{id}}                % identity map

 \newtheorem{theorem}{Theorem}[section]
 \newtheorem{conjecture}[theorem]{Conjecture}
 \newtheorem{proposition}[theorem]{Proposition}
 {
 \theorembodyfont{\rmfamily}
 \newtheorem{example}[theorem]{Example}
 \newtheorem{definition}[theorem]{Definition}
 \newtheorem{remark}[theorem]{Remark}
 }

\title{McKay correspondence}
\author{Miles Reid, Nagoya and Warwick}
\date{Feb 1997}
\begin{document}
\input epsf \epsfverbosetrue

\maketitle

 \section{Introduction}\label{sec:intro}

 \begin{conjecture}[since 1992]\label{conj:1992}
$G\subset\SL(n,\C)$ is a finite subgroup. Assume that the quotient $X=\C^n/G$
has a crepant resolution $f\colon Y\to X$ (this just means that $K_Y=0$, so
that $Y$ is a ``noncompact Calabi--Yau manifold''). Then there exist
``natural'' bijections
 \begin{align}
 \{\text{\em irreducible representations of $G$}\}
 &\to \text{\em basis of $H^*(Y,\Z)$} \\
 \{\text{\em conjugacy classes of $G$}\} &\to \text{\em basis of $H_*(Y,\Z)$}
 \end{align}

As a slogan ``$\text{representation theory of $G$}=\text{homology theory of
$Y$}$''.

Moreover, these bijections satisfy ``certain compatibilities''

 \begin{equation}
 \left.
 \begin{array}{rr}
\text{\em character table of $G$}\\
\text{\em McKay quiver}
 \end{array}
 \right\}
 \bij
 \left\{
 \begin{array}{ll}
\text{\em duality}\\
\text{\em cup product}
 \end{array}
 \right.
 \notag
 \end{equation}

\end{conjecture}

As you can see, the statement is still too vague because I don't say what
``natural'' means, and what ``compatibilities'' to expect. At present it
seems most useful to think of this statement as pointer towards the truth,
rather than the truth itself (compare Main Conjecture~\ref{conj:K}).

The conjecture is known for $n=2$ (Kleinian quotient singularities, Du Val
singularities). McKay's original treatment was mainly combinatorics
\cite{McK}. The other important proof is that of Gonzales-Sprinberg and
Verdier \cite{GSp-V}, who introduced the GSp--V or tautological sheaves,
also my main hope for the correspondence (1).

For $n=3$ a weak version of the correspondence (2) is proved in \cite{IR}.
We hope that a modification of this idea will work in general for (2); for
details, see \S\ref{sec:IR}.

 \paragraph{Contents} This is a rough write-up of my lecture at Kinosaki and
two lectures at RIMS workshops in Dec 1996, on work in progress that has not
yet reached any really worthwhile conclusion, but contains lots of fun
calculations. History of Vafa's formula, how McKay correspondence relates to
mirror symmetry. The main aim is to give numerical examples of how the McKay
correspondences (1) and (2) must work, and to restate
Conjecture~\ref{conj:1992} as a {\em tautology}, like the cohomology or
K-theory of projective space $\proj^n$ (see Main Conjecture~\ref{conj:K}).
Introduction to Nakamura's results on the Hilbert scheme of $G$-clusters.

 \paragraph{Credits} Very recent results of I. Nakamura on $G$-Hilb, who sent
me a first draft of \cite{N3} and many helpful explanations. Joint work with
Y.~Ito. Moral support and invaluable suggestions of S. Mukai. Support
Sep--Nov 1996 by the British Council--Japanese Ministry of Education exchange
scheme, and from Dec 1996 by Nagoya Univ., Graduate School of Polymathematics.

 \subsection{History} Around 1986 Vafa and others defined the {\em stringy
Euler number} for a finite group $G$ acting on a manifold $M$:
 \begin{equation}
 \begin{aligned}
 e_{\text{string}}(M,G)&=\text{crazy formula (you'd better forget it!)}\\
 &=\sum_{H\subset G} e(X_H) \times \card\{\text{conjugacy classes in $H$}\}.
 \end{aligned}
 \tag{$*$}
 \end{equation}

Here $X=M/G$, and $X$ is stratified by stabiliser subgroups: for a subgroup
$H\subset G$, set
 \begin{align*}
 M_H &=\{Q\in M | \Stab_G Q=H\},\\
   X_H &=\pi(M_H)\\
&=\{P\in X | \text{for $Q\in\pi\1(P)$, $\Stab_G Q$ is conjugate to $H$}\}.
 \end{align*}
The sum in ($*$) runs over all subgroups $H$, and $e(X_H)$ is the ordinary
Euler number. The mathematical formulation ($*$) is due to Hirzebruch--H\"ofer
\cite{HH} and Roan \cite{Roan}. If $M=\C^n$ and $G\subset\GL(n,\C)$ only fixes
the origin, then the closure of each $X_H$ is contractible, so that only the
origin $\{0\}=X_G$ contributes to the sum in ($*$), and
 \begin{equation}
 e_{\text{string}}(\C^n,G)=\card\{\text{conjugacy classes in $G$}\}.
 \notag
 \end{equation}

At the same time, Vafa and others conjectured the following:

 \begin{conjecture}[``physicists' Euler number conjecture'']\label{conj:vaf}
In appropriate circumstances,
 \begin{equation}
 e_{\text{\em string}}(M,G)=\text{\em Euler number of minimal resolution of
$M/G$.}
 \notag
 \end{equation}
 \end{conjecture}

The context is string theory of $M=\text{CY 3-fold}$, and the $G$ action on
$M$ is Gorenstein, meaning that it fixes a global basis
$s\in\om_M=\Oh(K_M)\iso\Oh_M$ (dualising sheaf $\om_M=\Om^3_M$). In
particular, for any point $Q\in M$, the stabiliser subgroup is in $\SL(T_QM)$.

At that time, the physicists possibly didn't know that there was a generation
of algebraic geometers working on minimal models of 3-folds, and possibly
naively assumed that in their cases, there exists a unique minimal resolution
$Y\to X=M/G$, so that $e_{\text{string}}(M,G)=e(Y)$. A number of smart-alec
\hbox{3-folders} raised various instinctive objections, that a minimal model
may not exist, is usually not unique etc.

However, it turns out that the physicists were actually nearer the mark. One
of the points of these lectures is that, in flat contradiction to the official
3-fold ideology of the last 15 years, in many cases of interest, there {\em
is} a distinguished crepant resolution, namely Nakamura's $G$-Hilbert scheme.

My guess of the McKay correspondences follow on naturally from Vafa's
conjecture, by the following logic. If $M=\C^n$, then one sees easily that
for any reasonable resolution of singularities $Y\to X=\C^n/G$, the
cohomology is spanned by algebraic cycles, so that
 \begin{equation}
 e(Y)=\sum H^{p,p}=\card\{\text{algebraic cycles of $Y$}\}.
 \notag
 \end{equation}

It seems unlikely that we could prove the numerical concidence
 \begin{equation}
 e(Y)=\card\{\text{conjugacy classes of $G$}\}
 \notag
 \end{equation}
without setting up some kind of bijection between the two sets.
\cite{IR} does so for $G\subset\SL(3,\C)$.

 \subsection{Relation with mirror symmetry, applications}
Consider:
 \begin{enumerate}
 \renewcommand{\labelenumi}{(\alph{enumi})}
 \item the search for mirror pairs;
 \item Vafa's conjecture;
 \item conjectural McKay correspondence;
 \item speculative theory of equivariant mirror symmetry ($G$-mirror
symmetry).
 \end{enumerate}
Historically, (a) led to (b), (b) led to (c), and logically (c) implies (b).
I have long speculated that (c) is connected to (d), and maybe even that
it would eventually be proved in terms of (d). The point is that up to now,
the known proofs of the McKay correspondence (even in 2 dimensions) rely on
the explicit classification of the groups, plus quite detailed calculations,
and it would be very interesting to get more direct relations.

I suggest below in \S\ref{sec:taut} that the McKay correspondence can be
derived in tautological terms. If this works, it will have applications to
(d). Some trivial aspects of this are already contained in Candelas and
others' example of the mirror of the quintic 3-fold \cite{C}, where you could
take intermediate quotients in the $(\Z/5)^3$ Galois tower. My suggestion is
that $G$-mirror symmetry should relate pairs of CYs with group actions, and
include the character theory of finite groups as the zero dimensional case. I
guess you're supposed to add an analog of ``complexified K\"ahler
parameters'' to the conjugacy classes, and ``complex moduli'' to the
irreducible representations. Another application (more speculative, this one)
might be to wake up a few algebraists.

 \subsection{Conjecture~\ref{conj:1992}, (1) or (2), which is better?} I
initially proposed Conjecture~\ref{conj:1992} in 1992 in terms of irreducible
representations, an analog of the formulations of McKay and of \cite{GSp-V}. I
was persuaded by social pressure around the Trento conference and by my
coauthor Yukari Ito to switch to (2); its advantage is that the two sides are
naturally graded, and we could prove a theorem \cite{IR}. Batyrev and
Kontsevich and others have argued more recently that (2) is the more
fundamental statement. However, the version of correspondence (2) in
cohomology stated in \cite{IR} gives a $\Q$-basis only: the crepant divisors
do not base $H^2(Y,\Z)$ in general: fractional combinations of them turn up
as $c_1(\sL)$ for line bundles on $Y$ that are eigensheaves of the group
action, that is, GSp-V sheaves for 1-dimensional representations of $G$.

These lectures return to (1), passing via K-theory; in this context, the
natural structure on the right hand side of (1) is not the {\em grading} of
$H^*$, but the {\em filtration} of $K_0Y$. In fact, my thoughts on (2) in
general are, to be honest, in a bit of a mess at present (see \S\ref{sec:IR}
and \S\ref{sec:Kexs} below).

 \section{First examples}
These preliminary examples illustrate the following points:
 \begin{enumerate}
 \item To construct a resolution of a quotient singularity $\C^n/G$, and a
very ample linear system on it, rather than {\em invariant rational
functions}, it is more efficient to use {\em ratios of covariants}, that is,
ratios of functions in the same character space. This leads directly to the
Hilbert scheme as a natural candidate for a resolution.
 \item Functions in a given character space $\rho$ define a tautological sheaf
$\sF_\rho$ on the resolution $Y\to X$, and in simple examples, you easily cook
up combinations of Chern classes of the $\sF_\rho$ to base the cohomology of
$Y$.
 \end{enumerate}

I fix the following notation: $G\subset\GL(n,\C)$ is a finite subgroup,
$X=\C^n/G$ the quotient, and $Y\to X$ a crepant resolution (if it exists).
For a given cyclic (or Abelian) group, I choose eigencoordinates
$x_1,\dots,x_n$ or $x,y,z,\dots$ on $\C^n$. I write
$\frac1{r}(a_1,\dots,a_n)$ for the cyclic group $\Z/r$ action given by
$x_i\mapsto \ep^{a_i}x_i$, where $\ep=\exp(2\pi i/r)=\text{fixed primitive
$r$th root of 1}$. Other notation, for example the lattice
$L=\Z^n+\Z\cdot\frac1{r}(a_1,\dots,a_n)$ of weights, and the junior simplex
$\De_{\text{junior}}\subset L_\R$ are as in \cite{IR}.

 \begin{example}\label{ex:A_n} The
 \begin{figure}[thb]
 \centering\mbox{\epsfbox{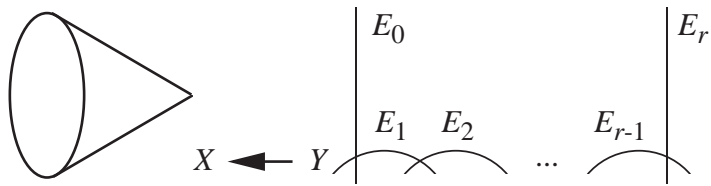}}
 \caption{$E_0$ and $E_r$ are the image of the $x$ and $y$ axes}
 \label{fig:A_n}
 \end{figure}
quotient singularity $\frac1{r}(1,-1)$. The notation means the cyclic
group $G=\Z/r$ acting on $\C^2$ by $(x,y)\mapsto(\ep x,\ep^{r-1}y)$. Everyone
knows the invariant monomials $u=x^r,v=xy,w=y^r$, the quotient map
 \begin{equation}
\C^2\to X=\C^2/G=\text{Du Val singularity $A_{r-1}$}:(uw=v^r)\subset\C^3,
 \notag
 \end{equation}
and the successive blowups that give the resolution $Y\to X$ and its chain of
$-2$-curves $E_1,\dots,E_{r-1}$ (Figure~\ref{fig:A_n}).
However, the new point to note is that each $E_i$ is naturally parametrised by
the ratio $x^i:y^{r-i}$. More precisely, an affine piece $Y_i\subset Y$ of the
resolution is given by $\C^2$ with parameters $\la,\mu$, and the equations
 \begin{equation}
x^i=\la y^{r-i},\quad y^{r-i+1}=\mu x^{i-1}\quad\text{and}\quad xy=\la\mu
 \label{eq:Ac}
 \end{equation}
define the $G$-invariant rational map $\C^2\broken Y_i$ (quotient map and
resolution at one go).

The ratio $x^i:y^{r-i}$ defines a linear system $|L(i)|$ on $Y$, with
intersection numbers
 \begin{equation}
L(i)\cdot E_j=\de_{ij}\quad\text{(Kronecker $\de$).}
 \notag
 \end{equation}
Thus, writing $\sL(i)$ for the corresponding sheaf or line bundle gives a
natural one-to-one correspondence from nontrivial characters of $G$ to
line bundles on $Y$ whose first Chern classes $c_1(\sL(i))\in H^2(Y,\Z)$
give the dual basis to the natural basis $[E_i]$ of $H_2(Y,\Z)$.
 \end{example}

 \begin{example}\label{ex:max}
One way of generalising Example~\ref{ex:A_n} to dimension 3. Let
 \begin{equation}
 G=\Span{\frac1{r}(1,-1,0),\frac1{r}(0,1,-1),\frac1{r}(-1,0,1,)}
 =(\Z/r)^2\subset\SL(3,\C)
 \notag
 \end{equation}
be the maximal diagonal Abelian group of exponent $r$. Then the first
quadrant of $L_\R$ has an obvious triangulation
 \begin{figure}[ht]
 \centering\mbox{\epsfbox{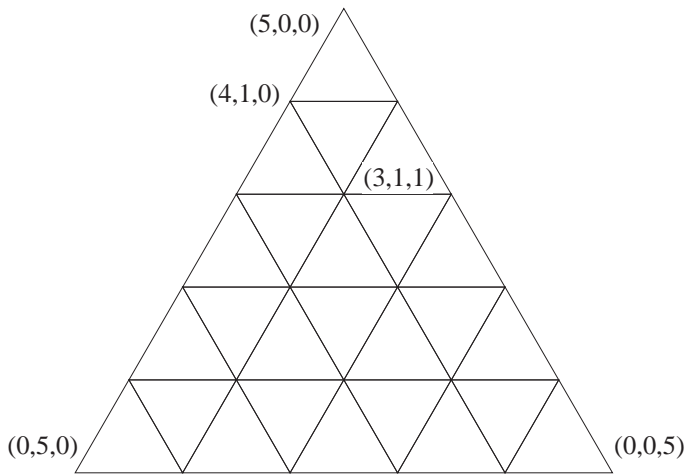}}
 \caption{Triangulation of $\protect\De_{\protect\text{junior}}$ in
Example~\ref{ex:max}}
 \label{fig:max1}
 \end{figure}
by regular simplicial cones that are basic for $L$ and have vertexes in the
junior simplex $\De_{\text{junior}}$. By toric geometry and the standard
discrepancy calculation \cite{YPG}, this triangulation defines a crepant
resolution $Y\to X=\C^3/G$.

From now on, restrict for simplicity to the case $r=5$ (featured on the
mirror of the quintic \cite{C}), whose triangulation is illustrated in
Figure~\ref{fig:max1}. $X=\C^3/G$ has lines of Du Val singularities
$A_4=\frac15(1,-1)$ along the 3 coordinate axes, the fixed locuses of the 3
generating subgroups $\frac15(1,-1,0)$ etc., of $G$.
 \begin{figure}[ht]
 \centering\mbox{\epsfbox{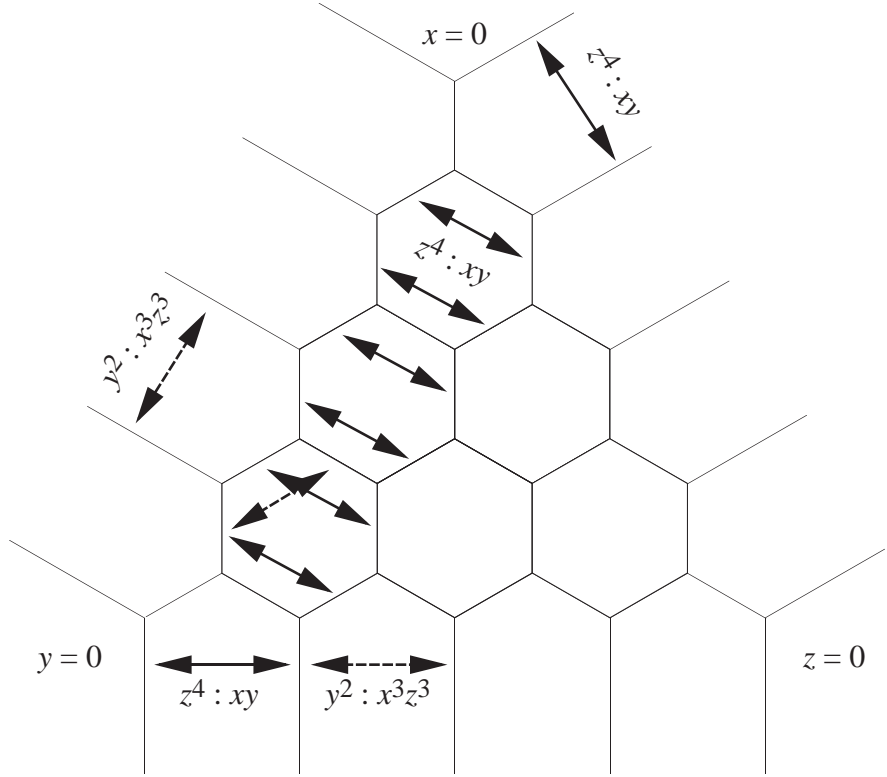}}
 \caption{The resolution corresponding to the triangulation of
Figure~\ref{fig:max1}}
 \label{fig:max2}
 \end{figure}
As illustrated in Figure~\ref{fig:max2}, the resolution $Y$ has 3 chains of 4
ruled surfaces over the coordinate axes of $X$, and 6 del Pezzo surfaces of
degree 6 (``regular hexagons'') over the origin. Every exceptional curve
stratum in the resolution is a $(-1,-1)$ curve.

Functions on the quotient $X=\C^3/G$ are given by $G$-invariant polynomials,
$k[X]=\C[x,y,z]^G$. To get more functions on $Y$ (and a projective embedding
of $Y$), consider the following ratios of monomials in the same eigenspace of
the $G$ action:
 \begin{equation}
  x^i:(yz)^{5-i} \quad\text{for $i=1,\dots,4$, and permutations of $x,y,z$.}
 \label{eq:ratio}
 \end{equation}
Each ratio (\ref{eq:ratio}) defines a free linear system on $Y$, and all
together, they define a relative embedding of $Y$ into a product of many
copies of $\proj^1$.
 \begin{figure}[hbt]
 \centering\mbox{\epsfbox{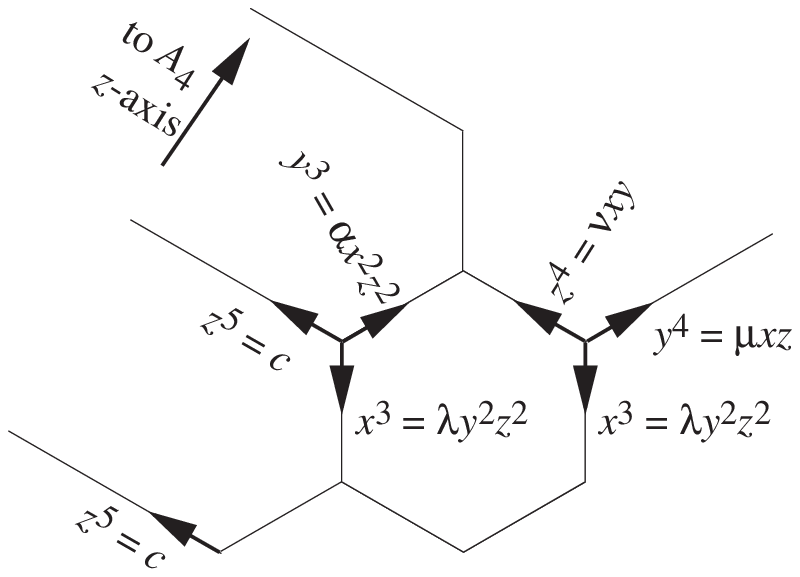}}
 \caption{Two affine pieces near the hexagon at (3,1,1)}
 \label{fig:max3}
 \end{figure}
For example, as shown in Figure~\ref{fig:max3}, the toric stratum at
$(2,2,1)$ is a del Pezzo surface of degree 6 embedded by the 3 ratios
$x^3:y^2z^2$, $y^3:x^2z^2$ and $z^4:xy$ (having product the trivial ratio
$1:1$). Figure~\ref{fig:max3} shows two affine pieces of $Y$, of which the
right-hand one is $\C^3$ with coordinates $\la,\mu,\nu$ related to $x,y,z$ by
a set of equations generalising (\ref{eq:Ac}):
 \begin{equation}
\begin{matrix}
x^3&=&\la y^2z^2\\
y^4&=&\mu xz\\
z^4&=&\nu xy
\end{matrix}
\qquad
\begin{matrix}
y^3z^3&=&\mu \nu x^2\\
x^2z^2&=&\la \nu y^3\\
x^2y^2&=&\la \mu z^3
\end{matrix}
\quad\text{and}\quad xyz=\la\mu\nu.
 \label{eq:maxr}
 \end{equation}

Denote the linear system $|x^i:(yz)^{5-i}|$ by $|L(x^i)|$, and similarly for
permutations of $x,y,z$. The sum of all the $|L(x^i)|$ is very ample on $Y$,
but their first Chern classes do not span $H^2(Y,\Z)$. To see this, recall
the del Pezzo surface $S_6$ of degree 6, the 3 point blowup of $\proj^2$
familiar from Cremona and Max Noether's elementary quadratic transformation.
It has 3 maps to $\proj^1$ and 2 maps to $\proj^2$; write
$e_1,e_2,e_3$ for the divisor classes of the maps to $\proj^1$, and $f_1,f_2$
for the maps to $\proj^2$. Then clearly,
 \begin{equation}
 \begin{gathered}
 e_1,e_2,e_3,f_1,f_2\quad \text{span} \quad H^2(S_6,\Z),\\
\text{with the single relation}\quad e_1+e_2+e_3=f_1+f_2.
 \end{gathered}
 \label{eq:reln}
 \end{equation}
For $S_6$ one of the hexagons of Figure~\ref{fig:max2}, the 3 maps to
$\proj^1$ are provided by certain of the linear systems $|L(x^i)|$. The two
maps to $\proj^2$ are provided by other character spaces: for example, for
the $(2,2,1)$ hexagon of Figure~\ref{fig:max3}, $f_1$ and $f_2$ are given by
the linear systems $|L(x^3y)|$ and $|L(xy^3)|$ corresponding respectively to
the ratios
 $$
\left(x^2z^4:x^3y:y^3z^2\right)
\quad\text{and}\quad
\left(xy^3:y^2z^4:x^3z^2\right)=
\left({1\over x^2z^4}:{1\over x^3y}:{1\over y^3z^2}
\right).
 $$
For each surface $S_6$, the generators $e_1,e_2,e_3,f_1,f_2$
correspond to certain characters of $G$. For example, if I choose the 3
generators
$\frac15(1,-1,0)$, $\frac15(0,1,-1)$ and $\frac15(-1,0,1)$ of $G$, the
characters of $x,y,z$ are
 $$
\renewcommand{\arraystretch}{1.2}
\begin{array}{rrr}
x & y & z
\\ \hline
1 & -1 & 0\\
0 & 1 & -1\\
-1 & 0 & 1
\end{array}
\quad\text{and my $(2,2,1)$ hexagon has} \quad
\begin{array}{ccc|cc}
e_1 & e_2 & e_3 & f_1 & f_2\\
x^3 & y^3 & z^4 & x^3y & xy^3\\ \hline
3 & 2 & 0 & 2 & 3\\
0 & 3 & 1 & 1 & 3\\
2 & 0 & 4 & 2 & 4
\end{array}
 $$
Moreover, you see easily that the relations (\ref{eq:reln}) actually hold
in $H^2(Y,\Z)$, not just in $H^2(S_6,\Z)$.

Represent each character of $G$ by a monomial $x^m$ (such as $x^i$ or
$x^3y$); this corresponds to a free linear system $|L(x^m)|$ on $Y$, in much
the same way as the $L(x^i:(yz)^{r-i})$ or $L(x^2z^4:x^3y:y^3z^2)$ just
described.

Now the McKay correspondence (1) of Conjecture~\ref{conj:1992} is the
following recipe:
 \begin{equation}
\text{monomial $x^m$} \mapsto \text{line bundle $\sL(x^m)$} \mapsto
c_1(\sL(x^m))\in H^2(Y,\Z).
 \notag
 \end{equation}
 These elements generate $H^2(Y,\Z)$, with one relation of the form
(\ref{eq:reln}) for every regular hexagon $S_6$ of the picture. Moreover,
each relation (\ref{eq:reln}) gives an element
 \begin{equation}
c_2(L(e_1)\oplus L(e_2)\oplus L(e_3))
-c_2(L(f_1)\oplus L(f_2))\in H^4(Y,\Z),
 \label{eq:reln2}
 \end{equation}
which is the dual element to $[S_6]\in H_4(Y,\Z)$. Indeed,
 \begin{align*}
c_2(L(e_1)\oplus L(e_2)\oplus L(e_3))\cdot S_6&=e_1e_2+e_1e_3+e_2e_3=3,\\
\text{and}\quad c_2(L(f_1)\oplus L(f_2))\cdot S_6&=f_1f_2=2.
 \end{align*}
I draw the McKay correspondence resulting from this cookery in
Figure~\ref{fig:max4}: each edge $E\iso\proj^1$ is labelled by the linear
system $L(x^m)$ with $L(x^m)\cdot E=1$, and each hexagon $S_6$ by 2 characters
corresponding to the two extra generators $f_1,f_2$ of $H^2(S_6,\Z)$ with the
relation which gives the dual element of $H^4(Y,\Z)$.
 \begin{figure}[ht]
 \centering\mbox{\epsfbox{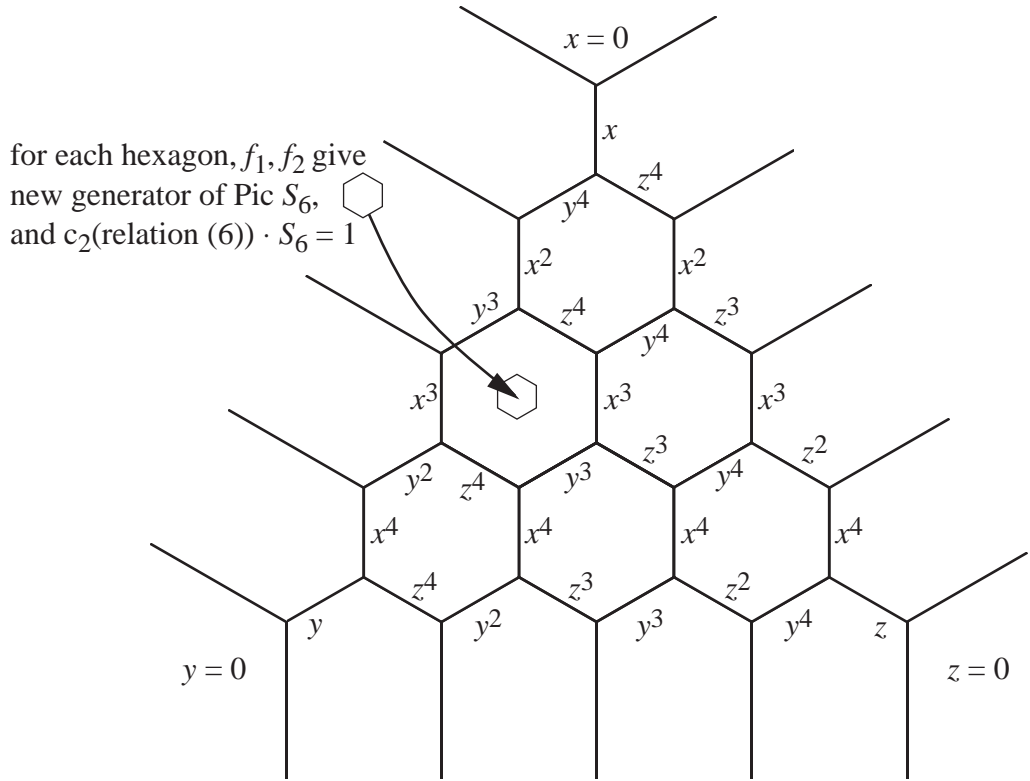}}
 \caption{McKay correspondence for Example~\ref{ex:max}}
 \label{fig:max4}
 \end{figure}
 \end{example}

One of the morals of this example is that we get a basis of cohomology in
terms of Chern classes of virtual sums of tautological bundles; this suggests
using the tautological bundles to base the K-theory of $Y$, and passing from
K-theory to cohomology by Chern classes or Chern characters. In fact, the
combinations used in (\ref{eq:reln2}) were fixed up to have zero first Chern
class, exactly what you must do if you want the second Chern character to
come out an integral class.

 \section{Ito--Reid, and the direct correspondence (2)}\label{sec:IR}
A group $G\subset\SL(n,\C)$ has a natural filtration by {\em age}. Namely, any
element $g\in G$ can be put in diagonal form by choosing $x_1,\dots,x_n$ to
be eigencoordinates of $g$. We write $g=\frac1{r}(a_1,\dots,a_n)$ to mean that
 \begin{equation}
g: (x_1,x_2,\dots,x_n) \mapsto (\ep^{a_1}x_1,\ep^{a_2}x_2,\dots,\ep^{a_n}x_n),
 \notag
 \end{equation}
where $\ep=\exp(2\pi i/r)=\text{fixed primitive $r$th root of 1}$, and
$a_i\in[0,1,\dots,n-1]$. Toric geometry tells us to consider the lattice
 \begin{equation}
 L=\Z^n+\Z \frac1{r}(a_1,\dots,a_n)
 \notag
 \end{equation}
(more generally for $A\subset G$ an Abelian group, we would add in lots of
vectors $\frac1{r}(a_1,\dots,a_n)$ for each $g\in A$). This consists of weightings
on the $x_i$, so that the invariant monomials have integral weights. Then for
any element $b=\frac1{r}(b_1,\dots,b_n)\in L$ with all $b_i\ge0$ (that is, $b$ in
the positive quadrant), define
 \begin{equation}
 \age(b)=\frac1{r}\sum b_i.
 \notag
 \end{equation}

In particular, for $g=\frac1{r}(a_1,\dots,a_n)$ in the unit cube,
 \begin{equation}
 \age(g)=\frac1{r} \sum a_i;
 \notag
 \end{equation}
this is obviously an integer (because $g\in\SL(n,\C))$ in the range $[0,n-1)$,
and this defines the age filtration.

Now any primitive vector $b=\frac1{r}(b_1,\dots,b_n)\in L$ and in the positive
quadrant defines a {\em monomial valuation} $v_b$ on the function field
$k(X)$ of $X$. Furthermore, the standard discrepancy calculation (see
\cite{YPG}) says that
 \begin{equation}
 \disc(v_b)=\age(b) - 1.
 \notag
 \end{equation}
\paragraph{Reminder:} The {\em discrepancy} $\disc v_b$ means that if I make a
blowup $W_b\to X$ so that $v_b$ is the valuation at a prime divisor
$F_b\subset W_b$, then $K_{W_b}=K_X+\disc(v_b) F_b$. Note also that {\em
junior} means $\age=1$, and {\em crepant} means $\text{discrepancy}=0$.
Any other questions?

The valuation $b$ defines a locus $E_b=\centre(v_b)\subset Y$. Consider only
weightings $b$ such that $v_b$ is the valuation of $E_b\subset Y$; this means
that if I blow up $Y$ along $E_b$, and $F_b$ is the exceptional divisor, then
$v_b$ is the valuation associated with the prime divisor $F_b\subset\widetilde
Y$. Since $Y$ is crepant, the adjunction formula for a blowup gives
 \begin{equation}
 \disc(v_b)=\codim E_b - 1,\quad\text{that is,}\quad\codim E_b=\age(b).
 \notag
 \end{equation}

In \cite{IR}, we uses this idea to give a bijection
 \begin{equation}
 \{\text{junior conjugacy classes of $G$}\} \to \{\text{crepant valuations of
$X$}\}
 \notag
 \end{equation}
which gave us a basis of $H^2(Y,\Q)$, and we dealt with $H^4(Y,\Q)$ by
Poincar\'e duality. Thus \cite{IR} only used the valuation theoretic
construction
 \begin{equation}
  b \mapsto v_b \mapsto E_b
 \notag
 \end{equation}
for $b$ in the junior simplex $\De_{\text{junior}}$. However, the same idea
obviously extends to give a correspondence from certain ``good'' elements $b$
to a set of locuses in $Y$ which generate $H_*(Y,\Z)$. Thus the idea for the
direct correspondence (2) is
 \begin{equation}
 \begin{aligned}
 G\ni g & \mapsto \text{collection of suitable $b$}\\
 & \mapsto \text{collection of locuses $E_b\subset Y$}.
 \end{aligned}
 \notag
 \end{equation}
The first step is by a mysterious cookery, which I only indicate by the
labelling in the two examples of \S\ref{sec:Kexs} below (it should be possible
to extract a good conjectural statement from this data).

 \section{Tautological sheaves and the main conjecture}\label{sec:taut} These
lectures are mainly concerned with providing experimental data for a suitably
rephrased Conjecture~\ref{conj:1992}, (1). In this section, I speculate on a
framework to explain what is going on, that might eventually lead to a proof.

The following is the main idea of \cite{GSp-V}. Given $G\subset\SL(n,\C)$, we
choose once and for all a complete set of irreducible representations
$\rho\colon G\to\GL(V_\rho)$. I use $\pi_*$ to view sheaves on $\C^n$ such as
the structure sheaf $\Oh_{\C^n}$ as sheaves on the quotient $\pi\colon\C^n\to
X$. Since $X$ is affine, these are really simply modules over
$k[X]=k[\C^n]^G$, so I usually omit $\pi_*$. Note that $k(\C^n)$ is a Galois
extension of $k(X)$, so that, by the cyclic element theorem of Galois theory,
it is the regular representation of $G$, that is, $k(\C^n)=k(X)[G]$; thus
$\pi_*\Oh_{\C^n}$ is generically isomorphic to the regular representation
$\Oh_X[G]$. For each
$\rho$, set
 \begin{equation}
 \sF'_\rho:=\Hom(V_\rho,\Oh_{\C^n})^G
 \notag
 \end{equation}
Then $\sF'_\rho\tensor V_\rho\subset\Oh_{\C^n}$ is the character subsheaf
corresponding to $V_\rho$; by the usual decomposition of the regular
representation, $\sF'_\rho$ is a sheaf of $\Oh_X$-modules of rank $\deg\rho$.
And there is a canonical decomposition
 \begin{equation}
 \Oh_{\C^n}=\sum_\rho \sF'_\rho\tensor V_\rho\quad\text{as $\Oh_X[G]$
modules.}
 \notag
 \end{equation}

Now let $f\colon Y\to X$ be a given resolution. Each $\sF'_\rho$ has a {\em
birational transform} $\sF_\rho$ on $Y$. This means that $\sF_\rho$ is the
torsion free sheaf of $\Oh_Y$ modules generated by $\sF'_\rho$, or if you
prefer, $\sF_\rho=f^*\sF'_\rho/(\text{torsion})$.

The sheaves $\sF_\rho$ are the {\em GSp-V sheaves}, or the {\em tautological
sheaves} of $Y$. Note that by definition, the $\sF_\rho$ are generated by
their $H^0$.

 \begin{conjecture}[Main conjecture]\label{conj:K} Under appropriate
circumstances, the\linebreak tautological sheaves $\sF_\rho$ form a $\Z$-basis
of the Grothendieck group $K_0(\Coh Y)$, and a certain cookery with their
Chern classes leads to a $\Z$-basis of $H^*(Y,\Z)$. A slightly stronger
conjecture is that the $\sF_\rho$ form a $\Z$-basis of the derived category
$D^b(\Coh Y)$.
\end{conjecture}

 \begin{remark} ``Appropriate circumstances'' in the conjecture include all
cases when $G\subset\SL(n,\C)$ and $Y=\GHilb$ is a crepant resolution. In
this case, these tautological sheaves $\sF_\rho$ have lots of good properties
(see \S\ref{sec:hilb}). But flops should not make too much difference to the
statement -- one expects a flopped variety $Y'$ to have more or less the same
homology and cohomology as $Y$, at least additively.
 \end{remark}

 \begin{example}\label{ex:Ver} $\frac1n(1,\dots,1)$ (with $n$ factors). The
quotient $X$ is the cone on the $n$th Veronese embedding of $\proj^{n-1}$, and
the resolution $Y$ is the anticanonical bundle of $\proj^{n-1}$, containing
the exceptional divisor $\proj^{n-1}$ with normal bundle
$\Oh(-n)=\om_{\proj^n}$. The tautological sheaves are
 \begin{equation}
 \Oh, \Oh(1),\dots,\Oh(n-1).
 \notag
 \end{equation}
That is, these are sheaves on $Y$ restricting down to the first $n$ multiples
of $\Oh(1)$ on $\proj^{n-1}$. It is well known that these sheaves form a
$\Z$-basis of the Grothendieck group $K_0(\proj^{n-1})$. It is a standard
(not quite trivial) bit of cookery with Chern classes and Chern characters to
go from this to a $\Z$-basis of $H^*(\proj^{n-1},\Z)$.
 \end{example}

 \begin{remark}
Recall the original (1977) {\em Beilinson diagonal trick}: the diagonal
$\De_{\proj^{n-1}}\subset\proj^{n-1}\times\proj^{n-1}$ is defined by the
section
 \begin{equation}
s_\De=\sum x_i'{\partial\ \over\partial x_i}\in
 p_1^*\Oh_{\proj^{n-1}}(1)\otimes p_2^*T_{\proj^{n-1}}(-1).
 \notag
 \end{equation}
Therefore, it follows (tautologically) that the derived category
$D^b(\Coh\proj^{n-1})$ (hence also the K theory $K_0$) has two ``dual'' bases
 \begin{equation}
 \Oh, \Om^1(1),\dots,\Om^{n-1}(n-1) \quad\text{vs.}\quad \Oh,
\Oh(-1),\dots,\Oh(-(n-1)).
 \notag
 \end{equation}
 \end{remark}

 \subsection*{Lame attempt to prove Conjecture~\ref{conj:K}}
 \paragraph{Step~I} The resolution $Y\to X$ is the quotient $A/H$ of an open
set $A\subset\C^N$ by a connected algebraic group $H$. In other words, by
adding extra variables in a suitable way, we can arrange to make the finite
quotient $X=\C^n/G$ equal to the quotient $\C^N/H$ of a bigger space by the
action of a connected group $H$ (the quotient singularities arise from jumps
in the stabiliser group of the $H$-action); moreover, we can arrange to
obtain the resolution $Y\to X$ by first deleting a set of ``unstable'' points
of $\C^N$ and then taking the new quotient $A/H$. For example, the Veronese
cone singularity of Example~\ref{ex:Ver} is $\C^{n+1}$ divided by
 \begin{equation}
\C^*\ni\la\colon(x_1,\dots,x_n;z)\mapsto(\la x_1,\dots,\la x_n;\la^{-n}z).
 \notag
 \end{equation}
(Obvious if you think about the ring of invariants). The finite group $\Z/n$
is the stabiliser group of a point of the $z$-axis. The blowup is the
quotient $A/\C^*$, where $A=\C^{n+1}\setminus\text{$z$-axis}$. (Because at
every point of $A$, at least one of the $x_i\ne0$, so the invariant ratios
$x_j/x_i$ are defined locally as functions on the quotient.)

 \paragraph{Step~II} Most optimistic form: the Beilinson diagonal trick may
apply to a quotient of the form obtained in Step~I. That is, the diagonal
$\De_Y\subset Y\times Y$ has ideal sheaf $\sI_{\De_Y}$ resolved by an exact
sequence in which all the other sheaves are of the form
$\sF_i\boxtimes\sG_i=p_1^*\sF_i\tensor p_2^*\sG_i$, where the $\sF_i$ and
$\sG_i$ are combinations of the tautological bundles.

It's easy enough to get an expression for the tangent sheaf of $Y$, in terms
of an Euler sequence arising by pushdown and taking invariants from the exact
sequence of vector bundles over $A$
 \begin{equation}
\Lie(H)\to T_A\to f^*(T_Y)\to 0,
 \label{eq:Eu}
 \end{equation}
where $\im\Lie(H)$ is the foliation by $H$-orbits. Maybe one can define a
filtration of this sequence corresponding to characters, and write the
equations of $\De_Y$ in terms of successive sections of twists of the graded
pieces. For example, the resolution $Y$ in Example~\ref{ex:Ver} is an affine
bundle over $\proj^{n-1}$, and the diagonal in $Y$ is defined by first taking
the pullback of the diagonal of $\proj^{n-1}$ (defined by the section
$\sum x_i'\partial/\partial x_i\in \Oh_{\proj^{n-1}}(1)\otimes
T_{\proj^{n-1}}(-1)$, the classic case of the Beilinson trick), then taking
the relative diagonal of the line bundle $\Oh(-n)$ over $\proj^{n-1}$.

 \paragraph{Step~III} The sheaves $\sF_i$ or $\sG_i$ appearing in a Beilinson
resolution form two sets of generators of the derived category $D^b(\Coh Y)$.
Indeed, for a sheaf on $Y$, taking $p_1^*$, tensoring with the diagonal
$\Oh_{\De_Y}$, then taking $p_{2*}$ is the identity operation. However, a
Beilinson resolution means that $\Oh_{\De_Y}$ is equal in the appropriate
derived category to a complex of sheaves of the form $\sF_i\boxtimes\sG_i$.
(This is a tautology, like saying that if $V$ is a vector space, and $f_i\in
V$, $g_i\in V^*$ elements such that $\id_V=\sum f_ig_i$, then $f_i$ and $g_i$
span $V$ and $V^*$.)

It should be possible to go from this to a basis of $D^b(\Coh Y)$ by an
argument involving Serre duality and the assumption $K_Y=0$. In this context,
it is relevant to note that the Beilinson trick leads to line bundles in the
range $K<\sF_i\le\Oh$ as one of the dual bases (for $\proj^{n-1}$, I believe
also in all the other known cases).

 \section{Generalities on $\protect\GHilb$}\label{sec:hilb}

The next sections follow Nakamura's ideas and results, to the effect that the
Hilbert scheme of $G$-orbits often provides a preferred resolution of
quotient singularities (see \cite{N1}--\cite{N3}, \cite{IN1}--\cite{IN3},
compare also \cite{N}, Theorem~4.1 and \cite{GK}); the results here are
mostly due to Nakamura. I write $M=\C^n$, and let $G\subset\GL(n,\C)$ be a
finite subgroup.

 \begin{definition}
$\GHilb $ is the fine moduli space of $G$-clusters $Z\subset M$.

Here a $G$-{\em cluster} means a subscheme $Z$ with defining ideal
$\sI_Z\subset\Oh_M$ and structure sheaf $\Oh_Z=\Oh_M/\sI_Z$, having the
properties:
 \begin{enumerate}
 \item $Z$ is a {\em cluster} (that is, a 0-dimensional subscheme). (Request
to 90\% of the audience: please suggest a reasonable translation of cluster
into Chinese characters (how about {\em tendan}, cf.\ {\em seidan} =
constellation, as in the Pleiades cluster?)
 \item $Z$ is $G$-invariant.
 \item $\deg Z=N=|G|$.
 \item $\Oh_Z\iso k[G]$ (the regular representation of $G$). For example, $Z$
could be a general orbit of $G$ consisting of $N$ distinct points.
 \end{enumerate}
 \end{definition}

 \begin{remark}
 \begin{enumerate}
 \item A quotient set $M/G$ is traditionally called an {\em orbit space}, and
that's exactly what $\GHilb M$ is -- the space of clusters of $M$ which are
scheme theoretic orbits of $G$.

 \item There is a canonical morphism $\GHilb M\to M/G$, part of the general
nonsense of Hilbert and Chow schemes: $\GHilb $ parametrises $Z$ by
considering the ideal $\sI_Z\subset\Oh_M$ as a point of the Grassmannian,
whereas the corresponding point of $M/G$ is constructed from the set of
hyperplanes (in some embedding $M\into\proj^{\text{large}}$) that intersect
$Z$.

 \item If $\pi\colon M\to M/G$ is the quotient morphism, and $P\in M/G$ a
ramification point, the scheme theoretic fibre $\pi^*P$ is always much too
fat; over such a point, a point of $\GHilb M$ adds the data of a subscheme $Z$
of the right length.

 \item I hope we don't need to know anything at all about $\Hilb^N M$ (all
clusters of degree $N=|G|$), which is pathological if $N,n\ge3$. Morally,
$\GHilb$ is a moduli space of points of $X=M/G$, and the right way to think
about it should be as a {\em birational change of GIT quotient} of $M/G$.

 \end{enumerate}
 \end{remark}

 \begin{conjecture}[Nakamura]\label{conj:N}
 \begin{enumerate}
 \renewcommand{\labelenumi}{(\roman{enumi})} 
 \item $Hilb^G M$ is irreducible.
 \item For $G\subset\SL(3,\C)$, $Y=\GHilb \C^3\to X=\C^3/G$ is a crepant
resolution of singularities. (This is mostly proved, see \cite{N3} and below.)
 \item For $G\subset\SL(n,\C)$, if a crepant resolution of\/ $\C^n/G$
exists, then $\GHilb\C^n$ is a crepant resolution.
 \item If $N$ is normal in $G$ and $T=G/N$ then $\Hilb^T\Hilb^N=\GHilb $.
 \end{enumerate}
 \end{conjecture}

 \begin{remark} For $n\ge4$, a crepant resolution $Y\to X$ usually does not
exist, but the cases when it does seem to be rather important. As Mukai
remarks, a famous theorem of Chevalley, Shephard and Todd says that for
$G\subset\GL(n,\C)$, the quotient $\C^n/G$ is nonsingular if and only if $G$
is generated by quasi\-reflections. Since we want to view $\GHilb \C^n$ as a
different way of constructing the quotient, the question of characterising
$G$ for which $\GHilb\C^n$ is nonsingular (or crepant over $\C^n/G$) is a
natural generalisation. We know that the answer is yes for groups
$G\subset\SL(2,\C)$, probably also $\SL(3,\C)$, so by analogy with
Shephard--Todd, I conjecture that it is also yes for groups generated by
subgroups in $G\subset\SL(2,\C)$ or $\SL(3,\C)$. For cyclic coprime groups
$\frac1{r}(a,b,c,d)$, based on not much evidence, I guess there is a crepant
resolution iff there are $\frac13(r-1)$ junior elements, that is, exactly one
third of the internal points of $\Box$ lie on the junior simplex (see
\cite{IR}); this is very rare -- by volume, you expect approx 4 middle-aged
elements for each junior one (as in most math departments). An easy example
to play with is $\frac1{r}(1,1,1,-3)$, which obviously has a crepant
resolution
 \begin{align*}
 \iff\enspace &\text{the simplex $\Span{([\frac r3],[\frac r3],[\frac r3],
r-3[\frac r3]), (1000),(0100),(0010)}$ is basic}\\
 \iff\enspace &r\equiv1\mod3.
 \end{align*}
For more examples, see also \cite{DHZ}.
 \end{remark}

 \begin{proposition}[Properties of $\protect\GHilb $]\label{prop:GHilb}
Assume Conjecture~\ref{conj:N}, (1). (In most cases of present interest, one
proves that $\GHilb$ is a nonsingular variety by direct calculation;
alternatively, if Conjecture~\ref{conj:N}, (1) fails, replace
$Hilb^G M$ by the irreducible component birational to $M/G$.)
 \begin{enumerate}
 \renewcommand{\labelenumi}{(\arabic{enumi})} 
 \item The tautological sheaves $\sF_\rho$ on $Y$ are generated by their
$H^0$.
 \item They are vector bundles.
 \item Their first Chern classes or determinant line bundles
 \begin{equation}
  \sL_\rho=\det \sF_\rho=c_1(\sF_\rho)
 \notag
 \end{equation}
define free linear systems $|L_\rho|$ according to (1), and are therefore nef.
 \item Any strictly positive combination $\sum a_\rho L_\rho$ of the $L_\rho$
is ample on Y.
 \item These properties characterise $\GHilb$ among varieties birational to
$X$ (or the irreducible component).
 \end{enumerate}
 \end{proposition}

 \begin{remark} If $G\subset\SL(n,\C)$ and $M=\C^n$, and $Y=\GHilb M$ is
nonsingular, the McKay correspondence says in particular that the $L_\rho$
span $\Pic Y=H^2(Y,Z)$ (this much is proved). In the 3-fold case, when $Y$ is
a crepant resolution, (3--4) resolve the contradiction with the expectation
of 3-folders, because they show how $\GHilb $ is distinguished among all
crepant resolutions of $X$. For if we flip $Y$ in some curve $C\subset Y$,
then by (4) we know that $LC>0$ for some $L=L_\rho$, and it follows that the
flipped curve $C'\subset Y'$ has $L'_\rho C'<0$. Thus (1--3) do not
hold on $Y'$.
 \end{remark}

\paragraph{Proof} Write $Y=\GHilb M$. By definition of the Hilbert scheme,
there exists a universal cluster $\sZ\subset Y\times M$, whose first
projection $p\colon\sZ\to Y$ is finite, with every fibre a $G$-cluster $Z$.
Now from the defining properties of clusters $p_*\Oh_Z$ is locally isomorphic
to $\Oh_Y[G]$, the regular representation of $G$ over $\Oh_Y$. In particular,
it is locally free, and therefore so are its irreducible factors
$\sF_\rho\otimes V_\rho$. Since $Z\subset M=\C^n$, the polynomial ring
$k[M]$ maps surjectively to every $\Oh_Z$, so that $p_*\Oh_Z$ is generated by
its $H^0$. This proves (1--3).

For any $G$-cluster $Z\in\GHilb M$, the defining exact
sequence
 \begin{equation}
 0\to\sI_Z\to\Oh_{C^n}\to\Oh_Z\to0
 \label{eq:defg}
 \end{equation}
splits as a direct sum of exact sequences (I omit $\pi_*$, remember):
 \begin{equation}
 0\to\sI_{Z,\rho}\to\sF'_\rho\tensor V_\rho\to F_{Z,\rho}\tensor V_\rho\to0
 \notag
 \end{equation}
Therefore $Z$ is uniquely determined by the set of surjective maps
$\sF_\rho\to F_{Z,\rho}$. This proves (4).

I now explain (5). The linear systems $|L_\rho|$ are birational in nature,
coming from linear systems of Weil divisors $|L_\rho|_X$ on the quotient
$X=M/G$, and their birational transforms on any partial resolution $Y'\to X$.
Now (5) says there is a unique model $Y$ on which these linear systems are
all free and their sum is very ample: namely, for a single linear system, the
blowup, and for several, the birational component of the fibre product of the
blowups. This also gives a plausibility argument for Conjecture~\ref{conj:N},
(iii): if we believe in the existence of one crepant resolution $Y'$, and we
admit the doctrine of flops from Mori theory, we should be able to flop our
way from $Y'$ to another model $Y$ on which the $|L_\rho|_Y$ are all free
linear systems. (This is not a proof: a priori, if the $L_\rho$ are dependent
in $\Pic Y$, a flop that makes one nef might mess up the nefdom of another.
However, it seems that the dependences are quite restricted; compare the
discussion at the end of Example~\ref{ex:II}.) \QED

I go through these properties again in the Abelian case, which is fun in its
own right, and useful for the examples in \S\ref{sec:Kexs}. Then an
irreducible representation $\rho$ is an element of the dual group 
 \begin{equation}
\widehat G=\bigl\{\text{homomorphisms $a\colon G\to r$th roots of 1 in
$\C^*$}\bigr\},
 \notag
 \end{equation}
where $r$ is the exponent of $G$. I write $\Oh_X(a)$ for the eigensheaf, and
$\sL_Y(a)$ for the tautological line bundle on $Y$ (previously $\sF'_\rho$
and $\sF_\rho$ respectively).

For any $Z$, the sequence (\ref{eq:defg}) splits as
 \begin{equation}
 0\to\bigoplus m_a \to\bigoplus \Oh_X(a) \to\bigoplus k_a \to 0
\quad\text{(sum over $a\in\widehat G$)},
 \notag
 \end{equation}
where $k_a$ is the 1-dimensional representation corresponding to $a$
(because of the assumption $\Oh_Z=k[G]$). Thus a $G$-cluster is exactly
the same thing as a set of maximal subsheaves
 \begin{equation}
 m_a\subset\Oh_X(a),\quad\text{one for every $a\in\widehat G$,}
 \notag
 \end{equation}
subject to the condition that $\sum m_a$ is an ideal in $\Oh_{C^n}$, that is,
that $m_a\Oh_X(b)\subset\Oh_X(a+b)$ for every $a,b\in\widehat G$.

Now it is an easy exercise to see that the Hilbert scheme parametrising
maximal subsheaves of $\Oh_X(a)$ is the blowup of $X$ in $\Oh_X(a)$, which I
write $\Bl_a X\to X$, and in particular, it is birational. It follows that
$\GHilb $ is contained in the product of these blowups:
 \begin{equation}
 \GHilb \subset\prod \Bl_a X
 \tag{$*$}
 \end{equation}
(where the product is the fibre product over $X$ of all the $\Bl_a X$ for
$a\in\widehat G$), and is the locus defined in this product by the ideal
condition:
 \begin{equation}
 m_a\Oh_X(b) \subset \Oh_X(a+b) \quad\text{for every $a,b\in\widehat G$}
 \tag{$**$}
 \end{equation}
(this obviously defines an ideal of $\Bl_a\times_X\Bl_b$).

By contruction of a blowup, each $\Bl_a$ has a tautological sheaf $\Oh_a(1)$,
which is relatively ample on $\Bl_a$. The tautological sheaves on $\GHilb $
are simply the restrictions of the $\Oh_a(1)$ to the subvariety ($*$). This
proves (1--4) again. \QED

 \begin{remark} The fibre product in ($*$) is usually reducible, with big
components over the origin (the product of the exceptional locuses of the
$\Bl_a$). However, it is fairly plausible that the relations ($**$) define an
irreducible subvariety. This is the reason for Conjecture~\ref{conj:N}, (1).
\end{remark}

 \section{Examples of Hilbert schemes}\label{sec:Kexs}
More experimental data, to support the following conclusions:
 \begin{enumerate}
 \renewcommand{\labelenumi}{(\alph{enumi})}
 \item $Y=\GHilb$ can be calculated directly from the definition; for 3-fold
Gorenstein quotients, it gives a crepant resolution, distinguished from other
models as embedded in projective space by ratios of functions in the same
character spaces.
 \item Conjecture~\ref{conj:1992} can be verified in detail in numerically
complicated cases. It amounts to a funny labelling by $a\in\widehat G$ of
curves and surfaces on the resolution.
 \item The relations in $\Pic Y$ between the tautological line bundles, whose
$c_2$ give higher dimensional cohomology classes, come from equalities
between products of monomial ideals.
 \end{enumerate}

 \begin{example}\label{ex:hirz} Examples~\ref{ex:A_n}--\ref{ex:max} are
$G$-Hilbert schemes. In fact the equations (\ref{eq:Ac}) and (\ref{eq:maxr})
were written out to define $G$-clusters.

Next, it is a pleasant surprise to note that the famous Jung--Hirzebruch
continued fraction resolution of the surface cyclic quotient singularity
$\frac1{r}(1,q)$ is the $G$-Hilbert scheme $(\Z/r)\text{-}\Hilb\C^2$. To save
notation, and to leave the reader a delightful exercise, I only do the example
$\frac15(1,2)$, where $5/2=[3,2]=3-1/2$; the invariant monomials and
weightings are as in Figure~\ref{fig:hirz}. As usual, $X=\C^2/G$ and $Y\to X$
is the minimal resolution, with two exceptional curves $E_1$ and $E_2$ with
$E_1^2=-2$, $E_2^2=-3$.
 \begin{figure}[ht]
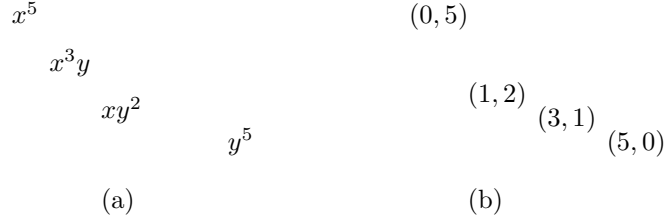

 $$
 \renewcommand{\arraycolsep}{2pt}
 \begin{array}{lllll}
 x^5\\[6pt]
 & x^3y \\[6pt]
 && xy^2 \\
 &&&\kern1cm y^5 \\[10pt]
 &&\text{(a)}
 \end{array}
 \kern2cm
 \renewcommand{\arraycolsep}{0pt}
 \begin{array}{lllll}
 (0,5)\\[18pt]
 && (1,2) \\[-3pt]
 &&& \kern4pt(3,1) \\[-3pt]
 &&&& \kern4pt(5,0) \\[10pt]
 &&\text{(b)}
 \end{array}
 $$
 \caption{Newton polygons (a) of invariant monomials and (b) of weights}
 \label{fig:hirz}
 \end{figure}
In toric geometry, $E_1$ corresponds to $(3,1)$ (as a vertex of the Newton
polygon (b) in the lattice of weights, or a ray of the fan defining the
resolution $Y$); the parameter along $E_1\iso\proj^1$ is $x:y^3$. Similarly,
$E_2$ corresponds to $(1,2)$ and has parameter $x^2:y$. Exactly as in
Figure~\ref{fig:A_n} and (\ref{eq:Ac}), a neighbourhood $Y_1$ of the point
$E_1\cap E_2$ is $\C^2$ with parameters $\la,\mu$, and the rational map
$\C^2\broken Y_1$ is determined by equations analogous to (\ref{eq:Ac}):
 \begin{equation}
 x^2=\la y,\quad y^3=\mu x,\quad\text{and}\quad xy^2=\la\mu.
 \label{eq:52clus}
 \end{equation}
These equations define a $G$-cluster $Z$: for a basis of
$\Oh_Z=k[x,y]/((\ref{eq:52clus}))$ is given by $1,y,y^2,x,xy$. Every
$G$-cluster is given by these equations, or by one of the following other two
types: $x^5=\la',y=\mu'x^2$ or $x=\la''y^3,y^5=\mu''$; the 3 cases correspond
to the 3 affine pieces with coordinates $\la,\mu$, etc. covering $Y$. The
generic $G$-cluster is $G\cdot(a,b)$ with $a,b\ne0$; all the equations
 \begin{equation}
x^5=a^5,\enspace x^3y=a^3b,\enspace xy^2=ab^2,\enspace
y^5=b^5,\enspace bx^2=a^2y,\enspace ay^3=b^3x
 \notag
 \end{equation}
vanish on $G\cdot(a,b)$, and since $a,b\ne0$, generators of its ideal can be
chosen in lots of different ways from among these, including the 3 stated
forms.

The ratio $x:y^3$ along $E_1$ and $x^2:y$ along $E_2$ define free linear
systems $|L(1)|$, $|L(2)|$ on $Y$ corresponding to the two characters $1,2$ of
$G=\Z/5$, with
 \begin{equation}
 \begin{aligned}
&L(1)\cdot E_1=1\\
&L(1)\cdot E_2=0
 \end{aligned}
\quad\text{and}\quad
 \begin{aligned}
&L(2)\cdot E_1=0\\
&L(2)\cdot E_2=1
 \end{aligned}
 \notag
 \end{equation}
These two give a dual basis of $H^2(Y,\Z)$, a truncated McKay correspondence.

 \paragraph{Exercise--Problem} The case of general $\frac1{r}(1,q)$ can be
done likewise; see for example \cite{R}, p.~220 for the notation, and compare
also \cite{IN2}. Problem: I believe that the minimum resolution of the other
surface quotient singularities is also a $G$-Hilbert scheme. The best way of
proving this may not be to compute $\GHilb$ exhaustively. In the $\SL(2,\C)$
case, Ito and Nakamura get the result $K_Y=0$ automatically, because the
moduli space $\GHilb$ carries a symplectic form.
 \end{example}

 \subsection*{The toric treatment of $\GHilb$} From now on, I deal mainly
with isolated Gorenstein cyclic quotient 3-fold singularities
$\frac1{r}(a,b,c)$, where $a,b,c$ are coprime to $r$ and $a+b+c=r$. If $G$ is
Abelian diagonal, then $X$ is obviously toric; however, it turns out that so
is the $G$-Hilbert scheme. There are two proofs; the better proof is that due
to Nakamura, described in \S\ref{sec:nak}. I now give a garbled sketch of the
first proof: I claim that the $G$-Hilbert scheme $\GHilb\C^n=Y(\Si)$ is the
toric variety given by the fan $\Si$, the ``simultaneous dual Newton polygon''
of the eigensheaves $\Oh_X(a)$, defined thus:
 \begin{quote}
for every character $a\in\widehat G$, write $\Oh_X(a)$ for the eigenspace of
$a$, $L(a)$ for the set of monomial minimal generators of $\Oh_X(a)$, and
construct the Newton polyhedron $\operatorname{Newton}(L(a))$ in the space of
monomials. Then $\Si$ is the fan in the space of weights consisting of the
cones $\Span{A_1,\dots,A_k}$ where the $A_i$ are weights having a common
minimum in every $L(a)$. This means that the 1-skeleton $\Si^1$ consists of
weights $A$ which either support a wall (= $(n-1)$-dimensional face) of
$\operatorname{Newton}(L(a))$ for some $a$, or which support positive
dimensional faces of a number of $L(a_j)$ whose product is $n-1$ dimensional
(in other words, ratios between monomials in the various $L(a_j)$ which are
minima for $A$ generate a function field of dimension $n-1$). Then
$\Span{A_1,\dots,A_k}$ is a cone of $\Si$ if and only if $\{A_i\}$ is a
complete set of weights in $\Si^1$ having a common minimum in every $L(a)$;
and $\Span{A_1,\dots,A_k}$ has dimension $d$ if and only if the ratio between
these minima span an $(n-d)$ dimensional space.
 \end{quote}
This definition is algorithmic, but quite awkward to use in calculations: you
have to list the minimal generators in each character space, and figure out
where each weight $A_i$ takes its least values; when $n=3$, you soon note that
the key point is the ratios like $x^3y:z^5$ between two monomials on an edge
of the Newton boundary.

 \begin{figure}[th]
 \centering\mbox{\epsfbox{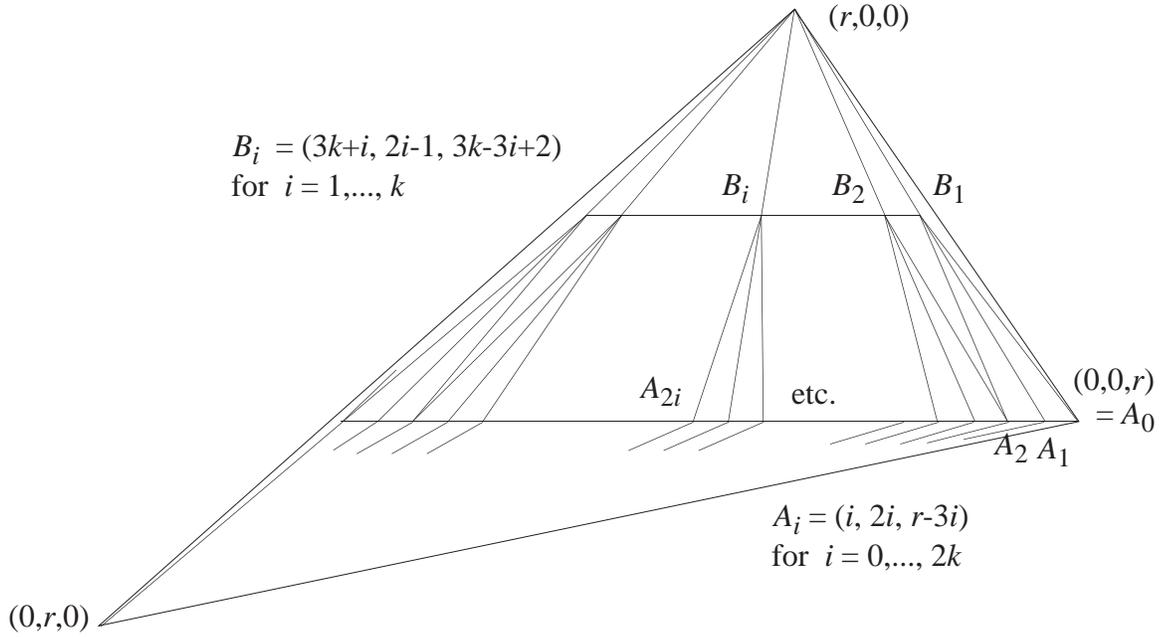}}
 \caption{$\protect\GHilb$ for $\protect\frac1{r}(1,2,-3)$. $B_i$ is joined to
$A_{2i-2},A_{2i-1},A_{2i}$}
 \label{fig:II}
 \end{figure}

 \paragraph{Sketch proof} Because $\Oh_Z=k[G]$ for $Z\in\GHilb$, for every
character $a$ of $G$, the generators of $L(a)$ map surjectively to the
1-dimensional character space $k_a$, so there is a well defined ratio between
the generators of $\sI_Z(a)$. This means that for fixed $Z$ and every $L(a)$,
we mark one monomial $s_a=x^{m(Z,a)}\in L(a)$ as the minimum of all the
valuations $A_1,\dots,A_k$ spanning a cone, and, using it as a generator, we
get the invariant ratios $x^{m'}/s_a$ as regular functions on $\GHilb$ near
$Z$.

 \begin{example}\label{ex:II} Consider $\frac1{r}(1,2,-3)$ where $r=6k+1$. The
quotient $X=\C^3/(\Z/r)$ is toric, and the $G$-Hilbert scheme is given by the
triangulation of the first quadrant of Figure~\ref{fig:II}.
 \begin{figure}[ht]
 \centering\mbox{\epsfbox{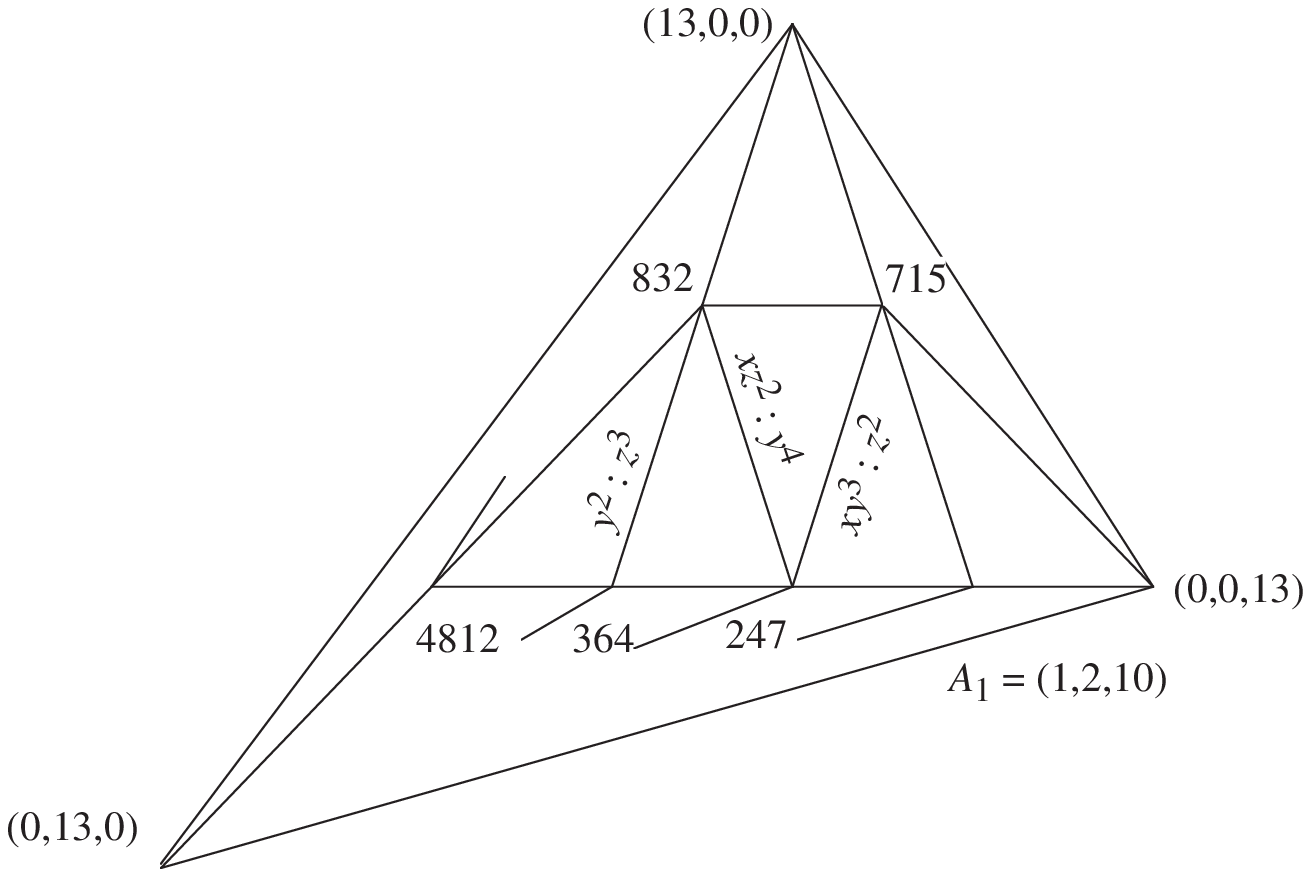}}
 \caption{$\protect\GHilb$ for $\protect\frac1{13}(1,2,10)$: why join
$(8,3,2)$---$(2,4,7)$?}
 \label{fig:13a}
 \end{figure}
This can be proved by carrying out the above proof explicitly. I omit the
laborious details, concentrating on one point: how does the Hilbert scheme
construction choose one triangulation in preference to another? For
simplicity, consider only $r=13$, so the triangulation simplifies to
Figure~\ref{fig:13a}. How do I know to join $(8,3,2)$---$(2,4,7)$ by a cone
$\si$, rather than $(7,1,5)$---$(3,6,4)$? By calculating $2\times2$ minors of
$\left(\begin{smallmatrix}8&3&2\\2&4&7\end{smallmatrix}\right)$, we see that
the parameter on the corresponding line $E_\si\in Y$ should be the ratio
$xz^2:y^4$, where $xz^2,y^4\in L(8)$. The Newton polygon of $L(8)$ is
 \begin{figure}[bht]
 \centering\mbox{\epsfbox{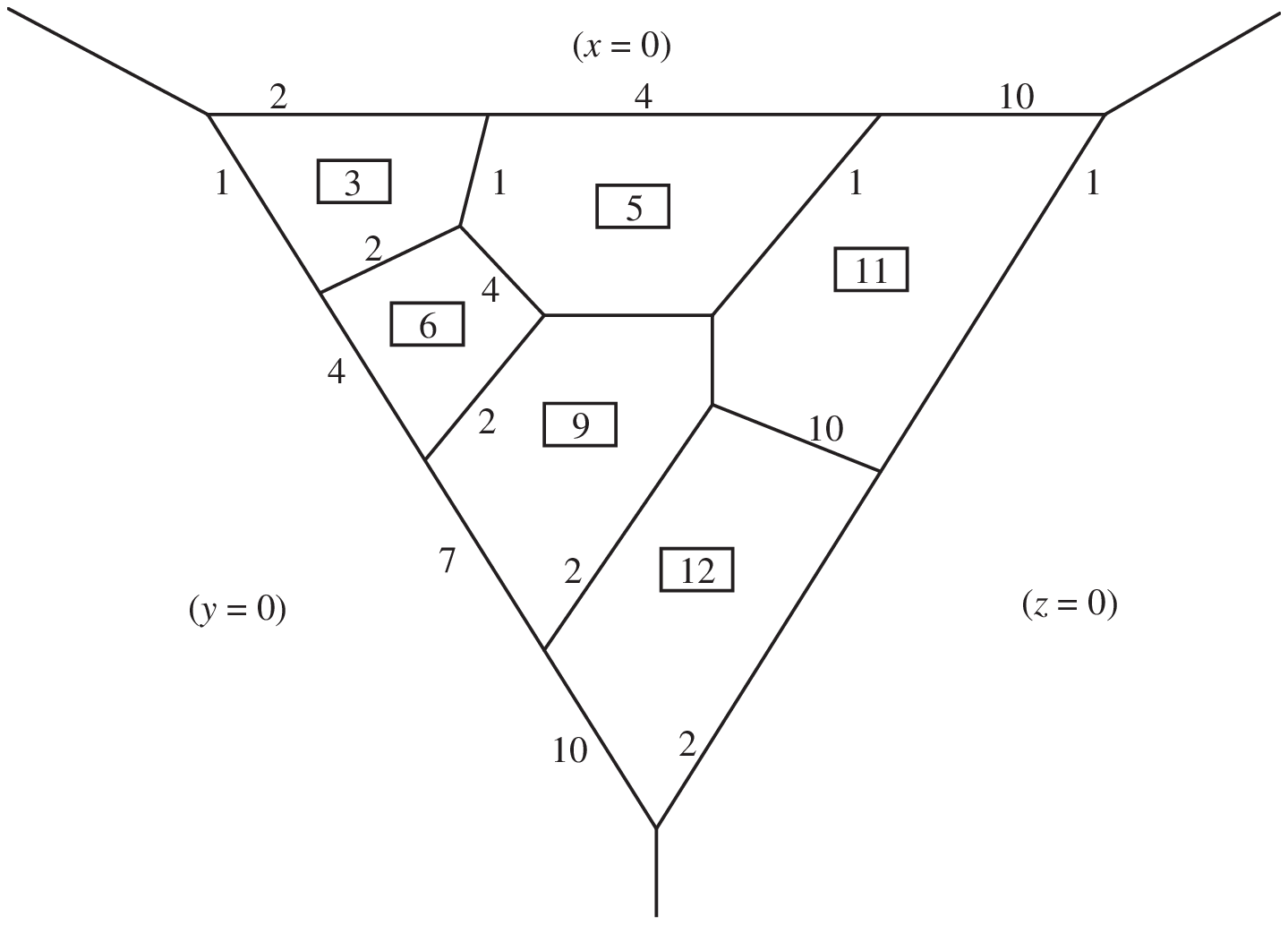}}
 \caption{The McKay correspondence for $\frac1{13}(1,2,10)$}
 \label{fig:13b}
 \end{figure}
 \end{example}
 $$
 \renewcommand{\arraystretch}{1.4}
 \begin{matrix}
x^8&x^6y&x^4y^2&x^2y^3&y^4\\
&&\kern-1cm(2,4,7)\kern-.6cm\\
&& xz^2&&y^2z^3\\
&&&\kern-1cm(8,3,2)\kern-.6cm\\
&&&& z^6
 \end{matrix}
 $$
(The figure is not planar: $xz^2$ and $y^4$ are ``lower''.) Here $(2,4,7)$
and $(8,3,2)$ have minima on the two planes as indicated, with common minima
on $xz^2$ and $y^4$, so that the linear system $|xz^2:y^4|$ can be free on
$L_\si$. But $(7,1,5)$ and $(3,6,4)$ don't have a common minimimum here:
$(7,1,5)$ prefers $y^4$ only, and $(3,6,4)$ prefers $xz^2$ only. If I join
$(7,1,5)$---$(3,6,4)$, the linear system $|xz^2:y^4|$ would have that line as
base locus.

The resolution is as in Figure~\ref{fig:13b}. The McKay correspondence marks
each exceptional stratum: a line $L$ parametrised by a ratio $x^{m_1}:x^{m_2}$
is marked by the common character space of $x^{m_1},x^{m_2}$. In other words,
a linear system such as $xz^2:y^4$ corresponds to a tautological line bundle
$\sL(xz^2:y^4)=\sL(8)$ with $c_1(\sL(8))\cdot L=1$.

The surfaces are marked by relations between the $c_1(\sL(i))$. In this case,
because there are no hexagons, these all arise from surjective maps
$\Oh_X(i)\otimes\Oh_X(j)\onto\Oh_X(i+j)$. For example, generators of the
character spaces $1,2,3$ are given by monomials (written out as Newton
polygons)
 \begin{equation}
\renewcommand{\arraycolsep}{3pt}
L(1):\enspace
 \begin{array}{lll}
 x & & y^7\\
 & y^2z\\
 z^4
 \end{array}, \quad
L(2):\enspace
 \begin{array}{llll}
 x^2 & & y\\
 xz^4\\
 z^8
\end{array} \quad\text{and}\quad
L(3):\enspace
 \begin{array}{llll}
 x^3 & xy & y^8\\
 x^2z^4 & y^3z\\
 xz^8\\
 z^{12}
\end{array} \quad
 \notag
 \end{equation}

 \begin{figure}[t!hb]
 \centering\mbox{\epsfbox{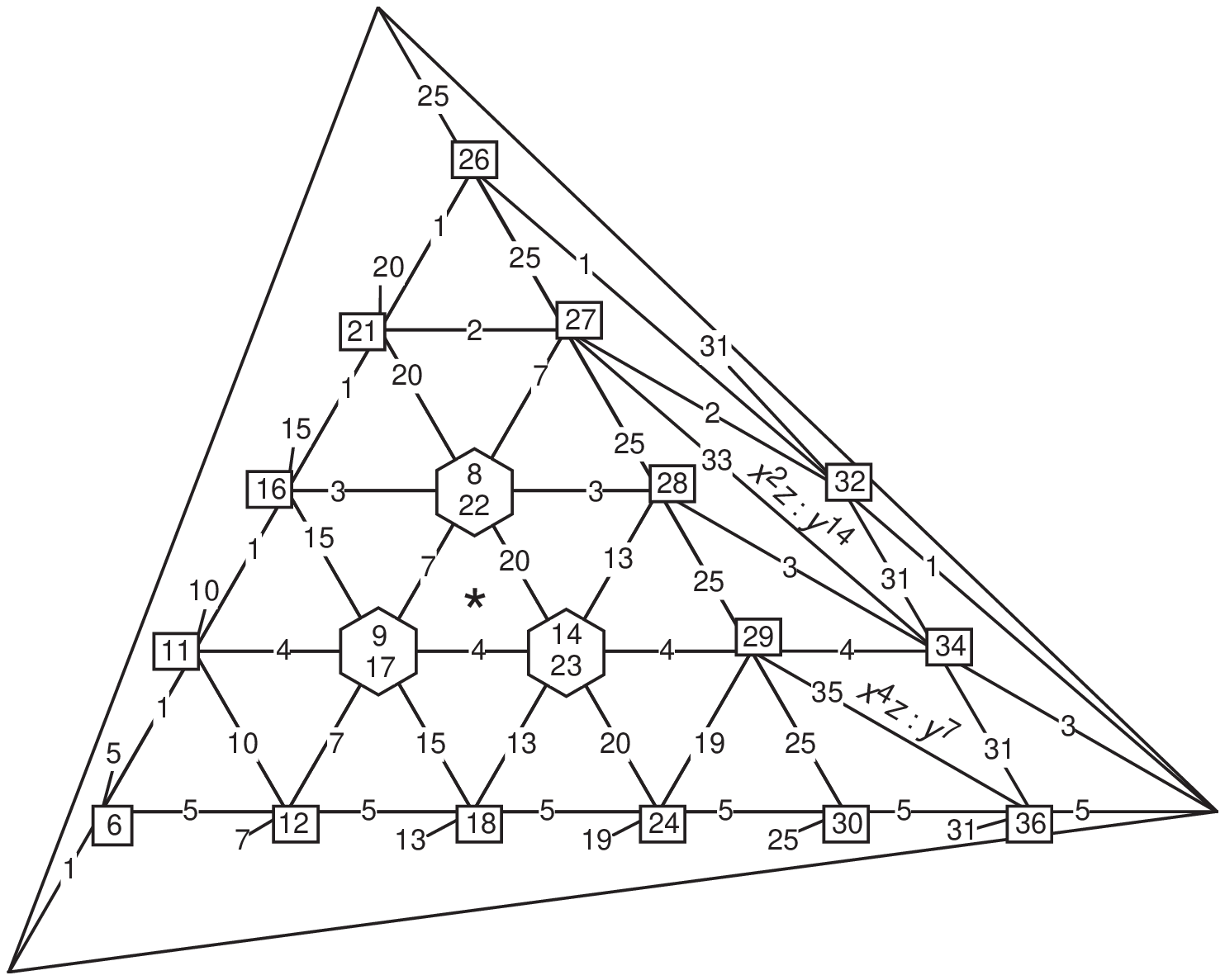}}
 \caption{The McKay correspondence for $\frac1{37}(1,5,31)$}
 \label{fig:37}
 \end{figure}

Clearly, $L(1)\otimes L(2)\onto L(3)$. (Thus this guy $\sL(3)$ is not active
in the resolution; in fact he's completely useless, so deserves to be
senior.) This means that on the resolution
 \begin{equation}
 c_1(\sL(3)-\sL(1)-\sL(2))=0,
 \notag
 \end{equation}
and $c_2(\sL(3)-\sL(1)-\sL(2))$ is the dual class to the top left surface in
Figure~\ref{fig:13b}.

 \begin{example}\label{ex:III} The $G$-Hilbert scheme for $\frac1{37}(1,5,31)$
is given by the triangulation in Figure~\ref{fig:37}, which also indicates
the labelling by characters of the McKay correspondence. I confine myself to
a few comments: on the right-hand side,
 \begin{align*}
 (1,5,31)\text{--}(9,8,20) \quad
 &\text{are joined by the ratio}\quad x^4z:y^7 \\
 (8,3,26)\text{--}(23,4,10) \quad
 &\text{are joined by the ratio}\quad x^2z:y^{14}
 \end{align*}
for reasons similar to those explained in Example~\ref{ex:II}. The
resolution has 3 regular hexagons (del Pezzo surfaces $S_6$), coming from the
regular triangular pattern on the left-hand side of Figure~\ref{fig:37}.
Tilings by regular hexagons appear quite often among the exceptional surfaces
of the Hilbert scheme resolution $Y$, as we saw in Figure~\ref{fig:max2}. The
reason for this is taken up again at the end of \S\ref{sec:nak}, see
Figure~\ref{fig:37b}. The cohomology classes dual to these 3 surfaces are
given as in (\ref{eq:reln2}) by taking $c_2$ of the relation
$e_1+e_2+e_3-f_1-f_2$, where the $f_1,f_2$ are the characters written in each
little hexagonal box of Figure~\ref{fig:37}, and $e_1,e_2,e_3$ are the
characters marking the 3 lines through the box. The relation
$e_1+e_2+e_3=f_1+f_2$ can also be expressed as equality between two products
of monomial ideals.
 \end{example}

 \section{Nakamura's proof that $\protect\GHilb$ is a crepant
resolution}\label{sec:nak}
 \begin{theorem}[Nakamura, very recent]\label{th:N}
For $G$ a finite diagonal subgroup of $\SL(3,\C)$, $Y=\GHilb\to X=\C^3/G$ is
a crepant resolution.
 \end{theorem}

\begin{figure}[h!b]
 \centering\mbox{\epsfbox{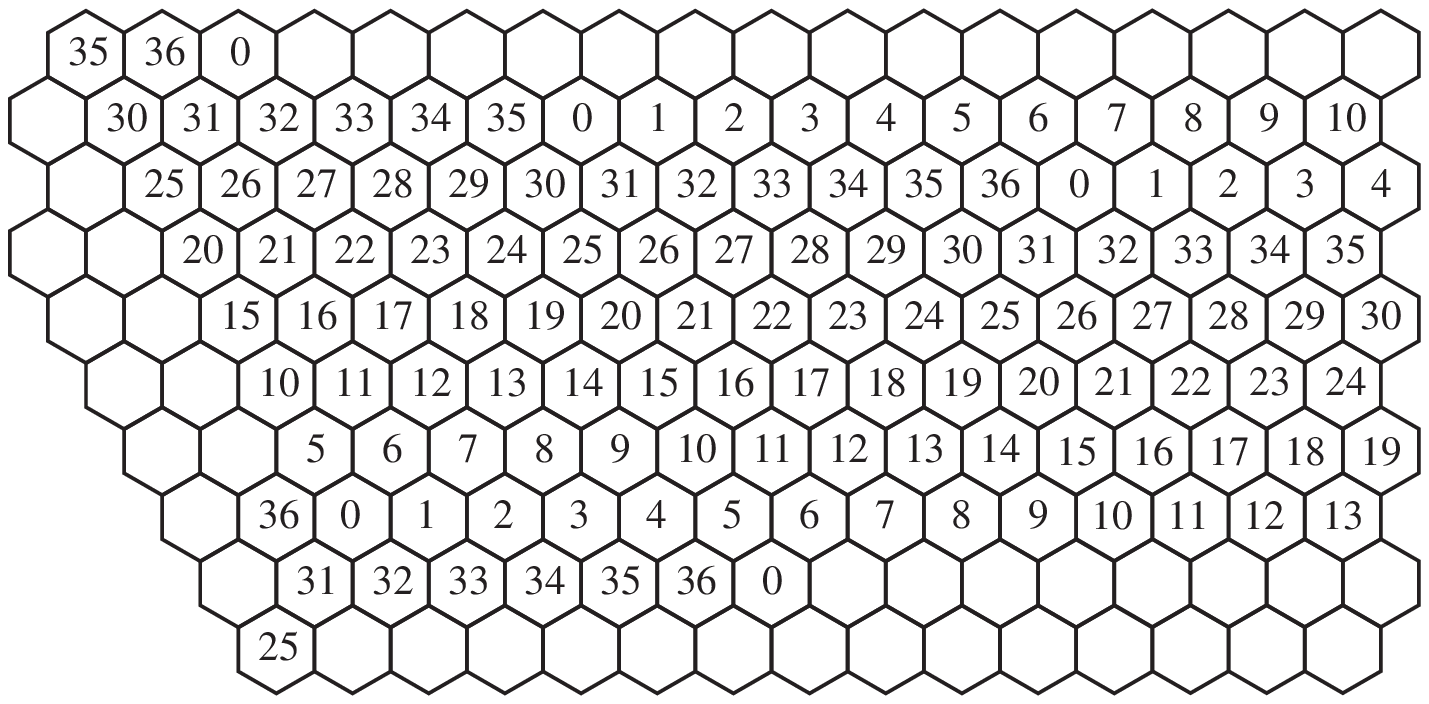}}
 \caption{The McKay quiver for $\frac1{37}(1,5,31)$.}
 \label{fig:honey37}
 \end{figure}

 \paragraph{Proof} I start from the {\em McKay quiver} of $G$ with the 3 given
characters $a,b,c$, corresponding to the eigencoordinates $x,y,z$, satisfying
$a+b+c=0$; to get the full symmetry, draw this as a doubly periodic
tesselation of the plane by regular hexagons, labelled by characters in
$\widehat G$:
 \begin{equation}
\renewcommand{\arraystretch}{1.3}
\renewcommand{\arraycolsep}{7pt}
\setcounter{MaxMatrixCols}{15}
\begin{matrix}
 & & & & & & & & & & & & & &\\[-18pt]
 &&& \kern-1cm \cdots \kern-1cm \\
 &&&& \kern-1cm 2b \kern-1cm \\
 &&&&& \kern-1cm b \kern-1cm && \kern-1cm a+b \kern-1cm 
 && \kern-1cm 2a+b \kern-1cm \\
\kern-1cm \cdots \kern-1cm &&
\kern-1cm 2b+2c \kern-1cm &&
\kern-1cm b+c \kern-1cm &&
\kern-1cm 0 \kern-1cm &&
\kern-1cm a \kern-1cm &&
\kern-1cm 2a \kern-1cm &&
\kern-1cm 3a \kern-1cm &&
\kern-1cm \cdots \kern-1cm \\
 &&&&& \kern-1cm c \kern-1cm && \kern-1cm a+c \kern-1cm \\
 &&&& \kern-1cm \cdots \kern-1cm
 \end{matrix}
 \label{fig:honey}
 \end{equation}
corresponding to the monomials
 \begin{equation}
\renewcommand{\arraystretch}{1.4}
\renewcommand{\arraycolsep}{8pt}
\setcounter{MaxMatrixCols}{13}
\begin{matrix}
 & & & & & & & & & & & &\\[-18pt]
 & \kern-1cm \cdots \kern-1cm \\
 && \kern-1cm y^2 \kern-1cm \\
 &&& \kern-1cm y \kern-1cm && \kern-1cm xy \kern-1cm &&
 \kern-1cm x^2y \kern-1cm \\
\kern-1cm y^2z^2 \kern-1cm &&
\kern-1cm yz \kern-1cm &&
\kern-1cm 1 \kern-1cm &&
\kern-1cm x \kern-1cm &&
\kern-1cm x^2 \kern-1cm &&
\kern-1cm x^3 \kern-1cm &&
\kern-1cm \cdots \kern-1cm \\
 &&& \kern-1cm z \kern-1cm && \kern-1cm xz \kern-1cm \\
 && \kern-1cm \cdots \kern-1cm
 \end{matrix}
 \notag
 \end{equation}
For $\frac1{37}(1,5,31)$, we get Figure~\ref{fig:honey37}; it is a quiver,
with arrows in the 3 principal directions ``add 1, 5 or 31''. Or you can view
it as the lattice of monomials modulo $xyz$, labelled with their characters
in the $\frac1{37}(1,5,31)$ action; then the arrows are multiplication by
$x,y,z$.

The whole of this business is contained one way or another in the hexagonal
figure (\ref{fig:honey}), together with its period lattice $\Pi$, and the many
different possible ways of choosing nice fundamental domains for the
periodicities; that is, we are doing Escher periodic jigsaw patterns on a
fixed honeycomb background. First of all, note that the periodicity of
(\ref{fig:honey}) is exactly the lattice of invariant Laurent monomials
modulo $xyz$. Call this $\Pi$.

The proof of Nakamura's theorem follows from the following proposition:

 \begin{proposition}\label{prop:N} For every $G$-cluster $Z$, the defining
equations (that is, the generators of $\sI_Z$) can be written as 7 equations
in one of the two following forms: either
 \begin{equation}
 \renewcommand{\arraycolsep}{1.5pt}
\begin{matrix}
x^{a+d+1}&=&\la y^bz^f\\
y^{b+e+1}&=&\mu z^cx^d\\
z^{c+f+1}&=&\nu x^ay^e
\end{matrix}
\qquad
\begin{matrix}
y^{b+1} z^{f+1}&=&\mu \nu x^{a+d}\\
z^{c+1} x^{d+1}&=&\la \nu y^{b+e}\\
x^{a+1} y^{e+1}&=&\la \mu z^{c+f}
\end{matrix}
\quad\text{and}\quad xyz=\la\mu\nu,\tag{$\uparrow$}
 \notag
 \end{equation}
for some $a,b,c,d,e,f\ge0$; or
 \begin{equation}
 \renewcommand{\arraycolsep}{1.5pt}
\begin{matrix}
x^{a+d}&=&\be\ga y^{b-1}z^{f-1}\\
y^{b+e}&=&\al\ga z^{c-1}x^{d-1}\\
z^{c+f}&=&\al\be x^{a-1}y^{e-1}
\end{matrix}
\qquad
\begin{matrix}
y^b z^f&=&\al x^{a+d-1}\\
z^c x^d&=&\be y^{b+e-1}\\
x^a y^e&=&\ga z^{c+f-1}
\end{matrix}
\quad\text{and}\quad xyz=\al\be\ga,\tag{$\downarrow$}
 \notag
 \end{equation}
for some $a,b,c,d,e,f\ge1$.
\end{proposition}

 \paragraph{Proof of Theorem~\ref{th:N}, assuming the proposition} Nakamura's
theorem follows easily, because $\GHilb$ is a union of copies of $\C^3$ with
coordinates $\la,\mu,\nu$ (or $\al,\be,\ga$), therefore nonsingular. Every
affine chart is birational to $X$, because it contains points with none of
$\la,\mu,\nu=0$ (or none of $\al,\be,\ga=0$). Moreover, an easy linear algebra
calculation shows that the equations ($\uparrow$) or ($\downarrow$) correspond to basic
triangles of the junior simplex, so that each affine chart of $\GHilb $ is
crepant over $X$. In more detail:

\paragraph{Case ($\uparrow$)} Write out the 3 x 3 matrix of exponents of the first
three equations of ($\uparrow$):
 $$
 \begin{matrix}
a+d+1&-b&-f\\
-d&b+e+1&-c\\
-a&-e&c+f+1
 \end{matrix}
 $$
(note that each of the 3 columns add to 1, more less equivalent to the junior
condition). The 2 x 2 minors of this give the 3 vertexes
 \begin{align*} P&=(bc+bf+ef+b+c+e+f+1,\quad ac+cd+df+d,\quad ab+ae+de+a),\\
Q&=(bc+bf+ef+b,\quad ac+cd+df+a+d+c+f+1,\quad ab+ae+de+e),\\
R&=(bc+bf+ef+f,\quad ac+cd+df+b,\quad ab+ae+de+a+b+d+e+1).
 \end{align*}

The triangle PQR ``points upwards'', in the sense that
 \begin{align*}
   &\text{$P$ is closest to $(1,0,0)$,}\\
   &\text{$Q$ is closest to $(0,1,0)$,}\\
   &\text{$R$ is closest to $(0,0,1)$.}
 \end{align*}
The 3 given ratios $x^{a+d+1}:y^bz^f$, etc.\ correspond to the 3 sides of
triangle $PQR$. In any case, all the vertexes belong to the junior simplex,
so that this piece of $\GHilb $ is crepant over $X$.

\paragraph{Case ($\downarrow$)} Write out the exponents of the second set of three
equations:
 $$
 \begin{matrix}
-(a+d)+1 &   b       &    f\\
    d    &  -(b+e)+1 &  c\\
    a    &      e    &  -(c+f)+1
 \end{matrix}
 $$
again, each of the 3 columns add to 1, and the $2\times2$ minors of this give
the 3 vertexes
 \begin{align*}
P&=(bc+bf+ef-b-c-e-f+1,     ac+cd+df-d,         ab+ae+de-a    ),\\
Q&=(    bc+bf+ef-b,     ac+cd+df-a-d-c-f+1,     ab+ae+de-e    ),\\
R&=(    bc+bf+ef-f,         ac+cd+df-b,     ab+ae+de-a-b-d-e+1),
 \end{align*}
all of which again belong to the junior simplex, so this affine chart is
also crepant over $X$. This time the triangle $PQR$ ``points downwards'', in
the sense that
\begin{align*}
  &\text{$P$ is furthest from $(1,0,0)$,}\\
  &\text{$Q$ is furthest from $(0,1,0)$,}\\
  &\text{$R$ is furthest from $(0,0,1)$.}
\end{align*}
The 3 given ratios $x^{a+d-1}:y^bz^f$, etc.\ again correspond to the 3 sides.
Q.E.D. for the theorem, assuming the proposition.

\paragraph{Proof of Proposition~\ref{prop:N}} Most of this is very geometric:
any reasonable choice of monomials in $x,y,z$ whose classes in $\Oh_Z$ form a
basis is given by a polygonal region $M$ of the honeycomb figure
(\ref{fig:honey}) satisfying 2 conditions:
 \begin{enumerate}
 \renewcommand{\labelenumi}{(\roman{enumi})} 
 \item in each of the 3 triants (triangular sector) it is concave, that is, a
downwards staircase: because it is a Newton polygon for an ideal;
 \item it is a fundamental domain of the periodicity lattice $\Pi$: because we
assume that $\Oh_Z=k[G]$, therefore every character appears exactly once.
 \end{enumerate}

The condition (ii) means that $M$ and its translates by $\Pi$ tesselate the
plane, so they form a kind of jigsaw pattern like the Escher periodic
patterns. However, in each of the 3 principal directions corresponding to the
$a$, $b$, and $c$-axes, there is only one acute angle, namely the summit at
the end of the $a$-axis (etc.). Therefore $M$ can only have one valley
(concave angle) in the $b,c$ triant. As a result, there is only one geometric
shape for the polygon $M$, the {\em tripod} or {\em mitsuya} (3 valleys, or 3
arrows) of Figure~I.

I introduce some terminology: the {\em tripod} $M$ has 3 {\em summits} at the
end of the axis of monomials $x^i$, and 3 {\em triants} or sectors of
$120^\circ$ containing monomials $x^iy^j$. Each triant has one {\em valley}
and two {\em shoulders} (incidentally, the 6 shoulders give the {\em socle} of
$\Oh_Z$).

 \begin{remark}\label{rem:deg} There are degenerate cases when some of the
valleys or summits are trivial (for example, $a=0$ in $\uparrow$). The most
degenerate case is a straight lines, when $\Oh_Z$ is based by powers of $x$
(say), and the equations boil down to $y=x^i,z=x^j$ (the $x$-{\em corner} of
the resolution). I omit discussion of these cases, since the equations of the
cluster $Z$ are always a lot simpler.
 \end{remark}

\begin{verbatim}
      I o o o
     o I o o o
    o o I o o o 
     o o I o o o o o o
      o o I o o o o o o
     o o o I I I I I I I     (Figure I)
    o o o I o o o o o o
     o o I o o
      o I o o
       I o o
\end{verbatim}

Thus there is only one ``geometric'' solution to the Escher jigsaw puzzle,
namely
\begin{verbatim}
               ...  I I I I I I I I I I
              u u u I u u u u u u u u u u u
               u u I u u
      I o o o   u I u u  
     o I o o o   I u u    I v v v 
    o o I o o o          v I v v v
     o o I o o o o o o  v v I v v v
      o o I o o o o o o  v v I v v v v v v ...   (Figure II)
     o o o I I I I I I I  v v I v
    o o o I o o o o o o  v v ...
     o o I o o          v v 
      o I o o
       I o o
\end{verbatim}

In particular, the external sides (going out to the 3 summits) are equal
plus-or-minus 1 to the opposite internal sides (going in to the 3 valleys).

 \begin{figure}[th]
 \centering\mbox{\epsfbox{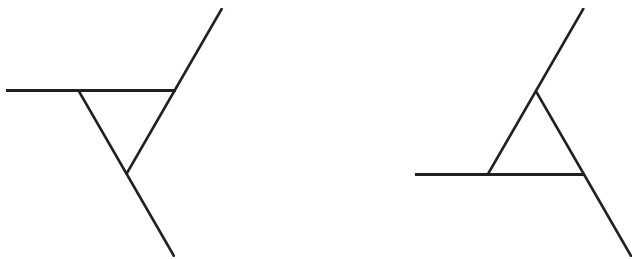}}
 \caption{Two different cocked hats}
 \label{fig:dmz}
 \end{figure}
However, the geometric statement of Figure~II is only exact for closed
polygons, whereas our tripods are Newton polygons spanned by integer points,
and are separated by a thin ``demilitarised zone'' between the integer points.
When you consider the tripods together with the integer lattices, there are
two completely different ways in which the three shoulders of neighbouring
tripods can fit together (corresponding to the two cyclic orders, or the two
cocked hats of Figure~\ref{fig:dmz}), namely either ($\uparrow$)
\begin{verbatim}
        y y y y z
         y y y z 
        x x x z z 
         x x x z z
\end{verbatim}
where the last $\tt y$ is just after the
last $\tt x$, and the shoulder of the
$\tt z$ is level with the top row of $\tt x$

or ($\downarrow$)
\begin{verbatim}
        y y y z z
         y y z z
        x x x z z
         x x x z z
\end{verbatim}
where the last $\tt y$ is just before the last $\tt x$ and the top row of
$\tt x$ is just below the shoulder of the $\tt z$.

The two different forms ($\uparrow$) and ($\downarrow$) come from this patching.

 \begin{remark}[Algorithm for $\GHilb$] Nakamura \cite{N3} gives an algorithm
to compute $\GHilb $ in this case as a toric variety. This can be viewed as a
way of classifying all the possible tripods in terms of elementary operations,
which correspond to the 0-strata and the 1-strata of the toric variety
$\GHilb $. You pass from an $\uparrow$ tripod to a $\downarrow$ one by
shaving off a layer of integer points one thick around one valley (assumed to
have thickness $\ge1$), and glueing it back around the opposite summit. And
vice versa to go from $\downarrow$ to $\uparrow$. You can start from
anywhere you like, for example from the $x$-corner (see Remark~\ref{rem:deg}).

Nakamura's algorithm applied to the statement in Proposition~\ref{prop:N}
expressed in terms of the fan triangulating the junior simplex, gives that if
$\uparrow$ and $a,b,c,d,\allowbreak e,f\ge k\ge2$ (say) then you can cross any
wall of the ``upwards'' triangle of the fan to get a new $\downarrow$
coordinate patch with $a',b',c',d',e',f'\ge k-1$, which corresponds to a
``downwards'' triangle, and vice-versa. It follows that the first triangle is
surrounded by a patch of width $k-1$ which is triangulated by the regular
triangular lattice, so that the resolution has a corresponding patch of
regular hexagons (that is, del Pezzo surfaces of degree 6).
Figure~\ref{fig:honey} shows the McKay quiver of $\frac1{37}(1,5,31)$ and
Figure~\ref{fig:37b} its fundamental domain giving the equations of
$G$-clusters
 \begin{equation}
 x^4=\la y^2z,\quad y^4=\mu xz^3,\quad z^5=\nu x^2y,\quad \text{etc.}
 \notag
 \end{equation}
on the coordinate chart of the resolution of $\frac1{37}(1,5,31)$,
corresponding to the starred triangle of Figure~\ref{fig:37}.

 \begin{figure}[ht]
 \centering\mbox{\epsfbox{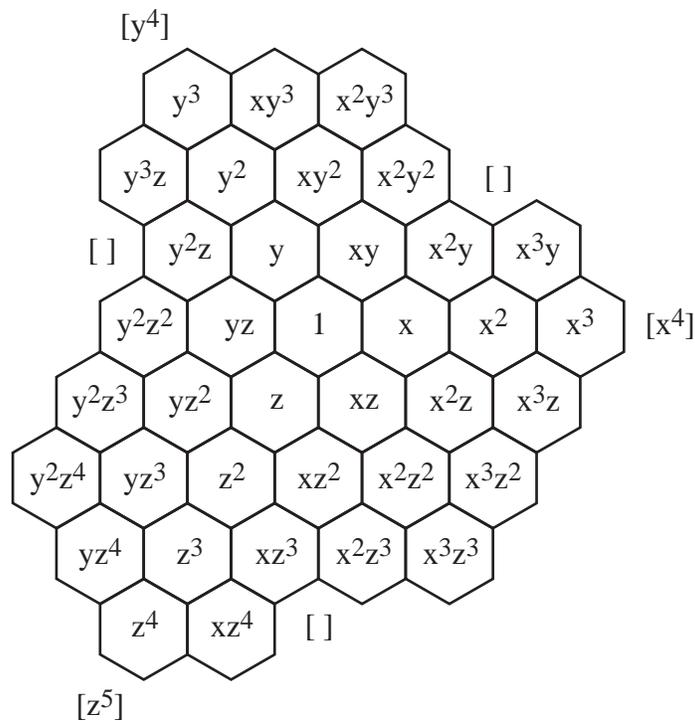}}
 \caption{A fundamental domain of the McKay quiver for $\frac1{37}(1,5,31)$}
 \label{fig:37b}
 \end{figure}

 \end{remark}

 \end{document}